\documentclass{thesis}

\usepackage{hyperref} % for clickable urls
\usepackage{graphicx}
\usepackage{amssymb} % for using \checkmark

\DeclareOldFontCommand{\sl}{\normalfont\slshape}{\@nomath\sl} % need this for compilation on arxiv.org, see https://www.tug.org/pipermail/tex-live/2016-July/038985.html

\newcommand{\thesistitle}{$J^+$-like Invariants under Bifurcations}
\newcommand{\thesisauthor}{Alexander Mai}
\newcommand{\submissiondate}{October 2022}

\title{\thesistitle}
\author{\thesisauthor}

\begin{document}

\setcounter{page}{1}

\pagestyle{intro}

{\Huge\textbf{\thesistitle}}

\begin{center}
	\textbf{\thesisauthor}, \href{mailto:alex.mai@posteo.net}{\Letter~alex.mai@posteo.net}\\
	\submissiondate, University of Augsburg\\
\end{center}

\section*{Abstract}
\addcontentsline{toc}{section}{\hspace{1.4em}Abstract}

We explore how the invariants~$J^+$, $J^-$, $\mathcal{J}_1$ and $\mathcal{J}_2$ of immersions -- generic (at most double points and only transverse intersections) planar smooth closed curves with non-vanishing derivative -- change under $k$-bifurcations ($k \ge 2 \in \mathbb{N}$), which are constructed by running $k$ times through an immersion and then perturbing it to be generic.

\section*{Introduction}
\addcontentsline{toc}{section}{\hspace{1.4em}Introduction}

An immersed loop is a smooth map $q: S^1 \to \mathbb{C} ,$ i.e. a smooth map of a circle into the plane. In this paper we identify the map with its image $K = q(S^1) \subset \mathbb{C} ,$ ignoring its parametrization and orientation. We call an immersed loop generic if it has only transverse self-intersections and all of them are double points. In this paper, all loops of interest are generic immersed loops, which we will simply call \emph{immersions} if not explicitly stated otherwise.

\emph{Vladimir Arnold} introduced three invariants for such immersions~\cite{arnold:paper}, of which the invariant $J^+$ is of special interest for this paper. If there are two immersions of the same rotation number, then by the \emph{Whitney--Graustein Theorem} one can be obtained from the other through a regular homotopy. Such two immersions have the same $J^+$ if and only if during a regular homotopy from one to the other, the number of positive direct self-tangencies is the same as the number of negative direct self-tangencies. A more thorough intro to the basics of the $J^+$-invariant and approaches to compute $J^+$ of interior sums can be found in my companion paper \emph{Introduction to Arnold's $J^+$-Invariant} \cite{mai:intro}.

We get a $k$-bifurcation $\widetilde{K}$ of some immersion $K$ if we run $k$ times through $K$ -- which results in a non-generic immersion as it intersects itself tangentially everywhere -- and then perturb it so that it is a generic immersion. Note that $k \ge 2, k \in \mathbb{N} .$ Depending on the way it is perturbed, there are arbitrary many different $k$-bifurcations $\widetilde{K}$ of $K .$

The goal of this paper is the following.

\begin{task}
	Find ways to calculate the $J^+$-invariant of $\widetilde{K}$, knowing the $J^+$-invariant of $K .$
\end{task}

To do this, we first observe that for any number $k$ there is no upper bound for the $J^+$-invariant due to slight variations that can create additional double points. But there is a lower bound:

\begin{theorema}
	$J^+(\widetilde{K}_k) \ge k^2 J^+(K) - (k^2 - k) .$
\end{theorema}

This lower bound is shown to be attained only by $k$-bifurcations with minimal number of double points, which are constructed and proven to be of minimal $J^+ .$ Our proof of Theorem~A has the following interesting consequence:

\begin{theoremb}
	$J^+(K) > 0 \quad \Rightarrow \quad J^+(\widetilde{K}_k) \ge k^2 + k .$
\end{theoremb}

Next we take into account how additional double points of a bifurcation are connected to $J^+$ and find, with $n_K$ the number of double points of an immersion $K$:

\begin{theoremc}
	$J^+(\widetilde{K}_k) = k^2 J^+(K) - (k^2 - 1) + n_{\widetilde{K}_k} - n_K k^2 .$
\end{theoremc}

Another one of the three invariants introduced by Arnold is the $J^-$-invariant. The $J^\pm$-invariants are closely related by the fact that their difference is the number of double points of the immersion, so for any immersion $K$ with $n_K$ double points we have:
$$ J^+(K) - J^-(K) = n_K . $$

For $J^-$ we see that the results simplify to a single formula independent of the number of double points $n_K .$ \\

The last one of the three invariants introduced by Arnold is called ``Strangeness'' -- abbreviated as~$St$ -- and of no further interest in this paper, as the discussed methods could not be applied to it. \\

Based on the $J^+$-invariant, \emph{Cieliebak, Frauenfelder and van Koert} introduce two new invariants~$\mathcal{J}_1$ and~$\mathcal{J}_2$ for any immersion $K \subset \mathbb{C}^* , \mathbb{C}^* \vcentcolon= \mathbb{C} \setminus \{ 0 \} ,$ see~\cite{kai:paper}. They are invariant under Stark--Zeeman homotopies (see~\cite{kai:paper} or Definition~\ref{def:szhom}) and rely on a chosen origin point, as they are dependent on the winding number~$\omega_0(K)$ of~$K$ around the origin:
$$ \mathcal{J}_1 = J^+(K) + \dfrac{\omega_0(K)^2}{2}, \quad \mathcal{J}_2 = \begin{cases}
	J^+(L^{-1} (K)), &\text{if \( \omega_0(K) \) odd} \\
	J^+(\widehat{K}), &\text{else,}
\end{cases}$$
with $L: \mathbb{C}^* \to \mathbb{C}^*, v \mapsto v^2$ the complex squaring map and $\widehat{K}$ one of the two up to rotation identical immersions of the preimage $L^{-1} (K) .$

Our results are shown to hold similarly for the invariants $\mathcal{J}_1$ and $\mathcal{J}_2 .$ Minor exceptions that apply for certain cases of the~$\mathcal{J}_2$-invariant are also shown. \\

This paper is organized as follows.
In Chapter~\ref{subsec:imbasics} we recall immersions and take note of some later useful observations on the winding number of components around a double point of an immersion.
Chapter~\ref{subsec:homevents} introduces the three disasters that are essential for Arnold's invariants.
Chapter~\ref{subsec:jplusminusinvs} is a quick overview over the properties of the $J^+$- and $J^-$-invariants and introduces \emph{Viro's formula}.
Chapter~\ref{subsec:szhom} is optional to any readers only interested in Arnold's original invariants, as we introduce the~$\mathcal{J}_1$- and~$\mathcal{J}_2$-invariants from \emph{Cieliebak, Frauenfelder and van Koert}.
Then we introduce bifurcations in Chapter~\ref{subsec:visbif} and discuss notable differences between bifurcations in Chapter~\ref{subsec:differentkbif}.
Finally all our results are presented from Chapter~\ref{subsec:kbifj} onwards, starting with results about the $J^+$-invariant. \\ %% TODO ALEX use this?

Acknowledgements: This paper is the result of my bachelor thesis, which would both not exist without the extraordinarily patient guidance and motivation from Kai Cieliebak, Urs Frauenfelder, Ingo Blechschmidt, Julius Natrup, Florian Schilberth, Milan Zerbin, Leonie Nießeler and other friends from the University of Augsburg.

\clearpage

\pagestyle{plain}

\clearpage

\section{Immersions and the invariants}
\label{sec:jinv}

In this chapter we will recall the necessary basics of immersions and homotopies and make some observations on winding numbers of connected components around double points and the concept of a double point's index. Then we introduce Arnold's $J^+$- and $J^-$invariant \cite{arnold:paper}, as well as the $\mathcal{J}_1$- and $\mathcal{J}_2$-invariants from \emph{Cieliebak, Frauenfelder and van Koert} \cite{kai:paper}.

The results around the~$J^+$-invariant do not use the other three invariants discussed in this paper. Therefore they can be ignored if one is only interested in~$J^+ ,$ or be saved for after the results around~$J^+$ have been digested. And while the results around these three invariants heavily use the results around~$J^+ ,$ they do not rely on each other. So for example Chapter~\ref{subsubsec:liftsdpsj2} on observations around the~$\mathcal{J}_2$-invariant is only needed for the results around~$\mathcal{J}_2$ -- there is no point in reading it for the results of~$J^-$ or~$\mathcal{J}_1 .$

A slightly more thorough introduction to $J^+$ on an undergraduate level can be found in my companion paper \emph{Introduction to Arnold's $J^+$-Invariant} \cite{mai:intro}.

\subsection{Immersions and winding numbers}
\label{subsec:imbasics}

\begin{defi}[Generic immersed loop]
	\label{def:imloop}
	An \emph{immersed loop} is a regular loop, i.e.\@ a smooth map
	$$q: S^1 \to \mathbb{C}$$
	with non-vanishing derivative, which we will, by slight abuse of notation, from here on identify with its oriented image $K = q(S^1) \subset \mathbb{C} .$
	We call it \emph{generic} if it only has \emph{transverse} self-intersections which are \emph{double points}.
\end{defi}

From here on we will often say \emph{immersion} to mean a \emph{generic immersed loop} if not stated otherwise.

We denote the winding number for all points in a connected component $C$ of $\mathbb{C} \setminus K$ as
$$\omega _C (K) \vcentcolon= \omega _y (K), \quad y \in C \text{ arbitrary} .$$

\begin{figure}[h!]
	\centering
	\includegraphics[scale=0.6]{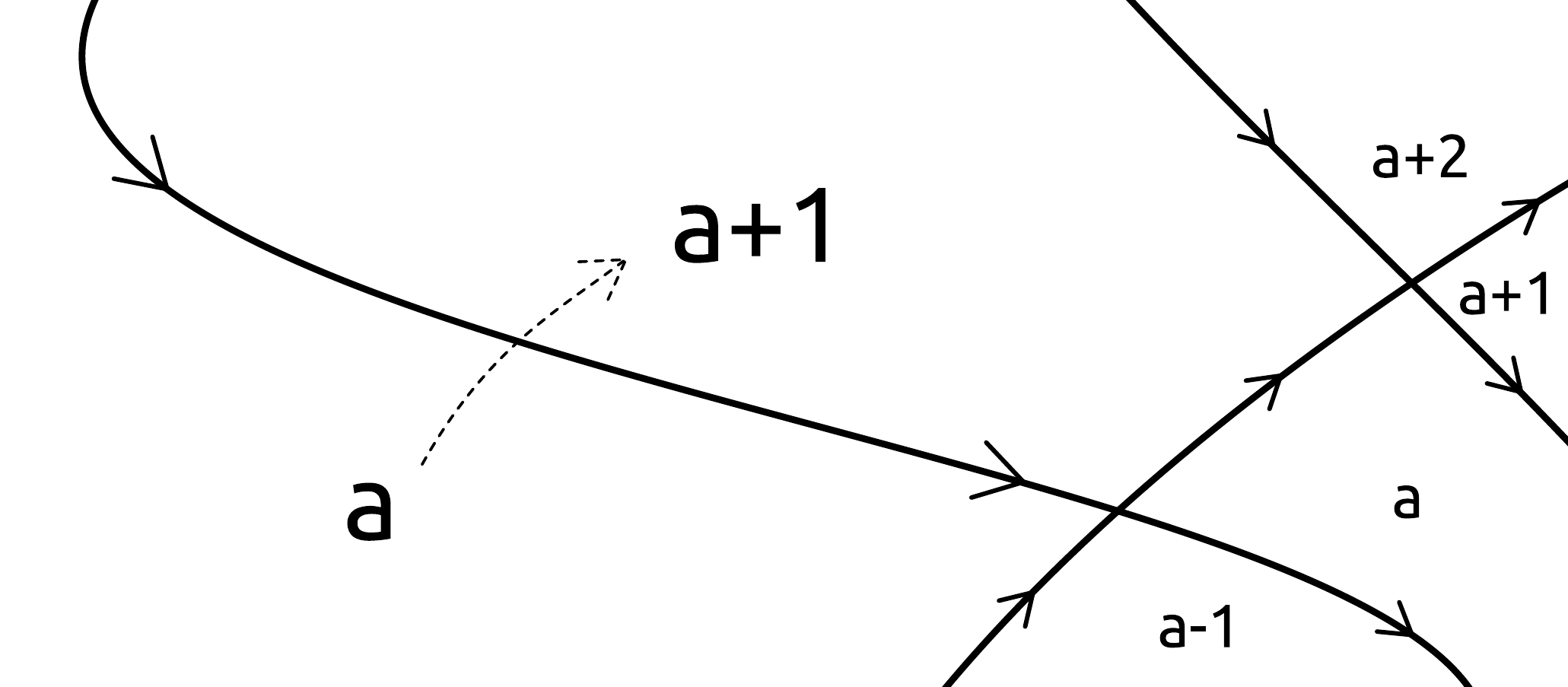}
	\caption{Changing winding number between adjacent connected components.}
	\label{fig:windingstep}
\end{figure}

The difference of the winding number of two connected components of $\mathbb{C} \setminus K$ that are adjacent to each other is always equal to $1.$ Looking at any connected component $C$ with a winding number of $a,$ the winding number of any adjacent component $C'$ is $a + 1$ ($a - 1$) if the immersion's arc, looking from $C$ to $C',$ is oriented from left to right (right to left). In other words, a flatlander living in one of the connected components can traverse an arc of the immersion to get to another connected component and the winding number of her location will increase (decrease) by $1$ if the traversed immersion arc is oriented from left to right (right to left), see Figure~\ref{fig:windingstep}.

With this trick it is easy to label all components' winding number in a picture of an immersion, starting from the unbounded component with winding number $0.$

\begin{thesislemma}
	\label{lem:dpwinding}
	The winding numbers around a double point $p$ are always of this form: let the lowest winding number around $p$ be $a \in \mathbb{Z}.$ Only one of the four components has winding number $a$, the component on the opposite site has winding number $a + 2$ and the other two have winding number~$a + 1 .$
\end{thesislemma}

\begin{figure}[h!]
	\centering
	\includegraphics[scale=0.45]{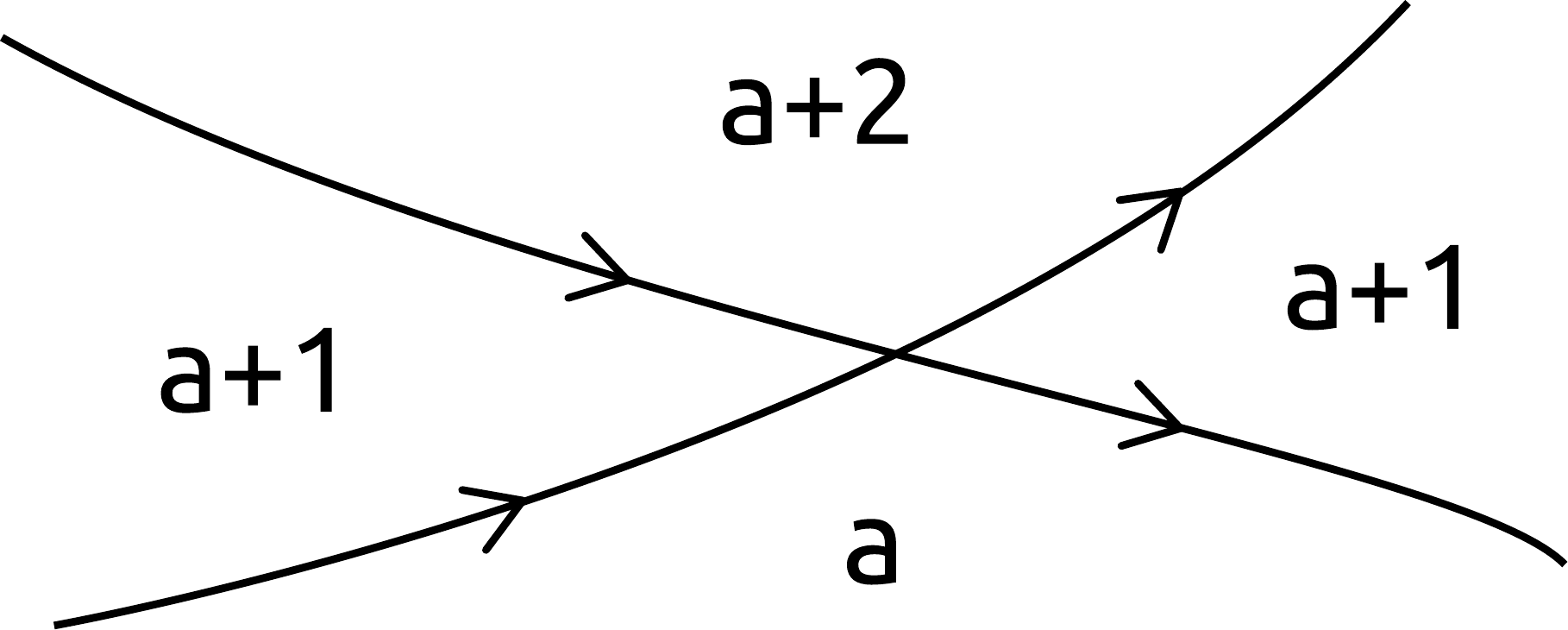}
	\caption{Winding numbers around a double point.}
	\label{fig:dpwinding}
\end{figure}

\begin{proof}
	See Figure~\ref{fig:dpwinding}. If $a$ is the lowest winding number around $p,$ the two adjacent components around $p$ have to be $a + 1,$ so the arcs to cross to get to these components have to be oriented left to right. This means the other two arcs also have to be oriented left to right and the last component has winding number $a + 2 .$
\end{proof}

\begin{defi}[Index $\operatorname{ind}_p(K)$]
	\label{def:indexdp}
	The \emph{index $\operatorname{ind}_p(K)$ of a double point $p$ in $K$} is the arithmetic mean of the winding numbers of the four connected components of $\mathbb{C} \setminus K$ adjacent to $p .$

	When calculating $\operatorname{ind}_p(K)$ and two corners of an adjacent component are equal to $p ,$ we count this one component and its winding number twice in the arithmetic mean to get a total of four adjacent components.
\end{defi}
\vspace*{0.4em}

\begin{remmm}
	\label{rem:indexdptwice}
	Using Lemma~\ref{lem:dpwinding} we can see that $\operatorname{ind}_p(K)$ is always equal to the winding number $a + 1$ appearing twice around $p,$ with $a \in \mathbb{Z} \,$ the lowest winding number around $p$ (same as in the lemma), as
	$$\dfrac{(a + (a + 1) + (a + 1) + (a + 2))}{4} = a + 1 .$$
\end{remmm}
\vspace*{0.4em}

\begin{remmm}
	\label{rem:indexdpside}
	It follows that whenever we draw a double point with the intersecting parts of the immersion both directed to the right (in the picture), then the index of that double point is equal to the winding number of the connected component to its right (or left). See Figure~\ref{fig:dpwinding}.

	This is also true if the intersecting parts are both directed from right to left.
\end{remmm}
\vspace*{0.4em}

\begin{defi}[Rotation number $\operatorname{rot}$]
	The \emph{rotation number} $\operatorname{rot}(q)$ of an immersed loop $q : S^1 \to \mathbb{C}$ is equal to the winding number of $\dot{q},$ its derivative, around the origin. So $\operatorname{rot}(q) = \omega_0(\dot{q}).$
\end{defi}

\subsection{Events during homotopies of immersed loops}
\label{subsec:homevents}

\begin{defi}[Regular homotopy]
	A \emph{regular homotopy} between two immersed loops $q$ and $q'$ is a smooth map $$h : S^1 \times [0, 1] \to \mathbb{C},$$ with $h(\cdot, 0) = q$ and $h(\cdot, 1) = q'$ and $h(\cdot, t) : S^1 \to \mathbb{C}$ an immersed loop $\forall t \in [0, 1]$.

	We call two immersed loops \emph{regularly homotopic} if there is a regular homotopy between these two immersed loops.
\end{defi}

By the \emph{Whitney--Graustein Theorem} any two immersions are regularly homotopic if and only if they have the same rotation number.

\begin{remark}
	The immersions during a regular homotopy between two generic immersions do not have to always be generic. There can be $t_0 \in (0, 1)$ with $h(\cdot, t_0)$ not generic, i.e.\@ with self-tangencies or at least triple points. Without loss of generality, these moments $t_0$ are isolated within our regular homotopies and even in the non-generic case we have at most triple points. This is honest, because if the immersion during a homotopy is not generic in all but isolated moments, small perturbations of the homotopy can make the immersion generic in all but isolated moments.
\end{remark}

There are three important isolated scenarios that can happen an arbitrary amount of times during a regular homotopy -- with only isolated double points for tangential intersections and only isolated triple points for transverse intersections.

\vspace*{0.2em}
\begin{defi}[Direct (inverse) positive (negative) self-tangency]
	\label{def:selftang}
	A \emph{self-tangency} is the event of an immersion crossing itself tangentially during a regular homotopy. We call this self-tangency \emph{direct} (\emph{inverse}) if both parts involved in the crossing of the immersion are (are not) oriented in the same direction.

	We call a direct self-tangency \emph{positive} (\emph{negative}) if the number of double points increases (decreases) by $2 .$
\end{defi}
\vspace*{0.8em}

\begin{figure}[h!]
	\centering
	\includegraphics[scale=0.3]{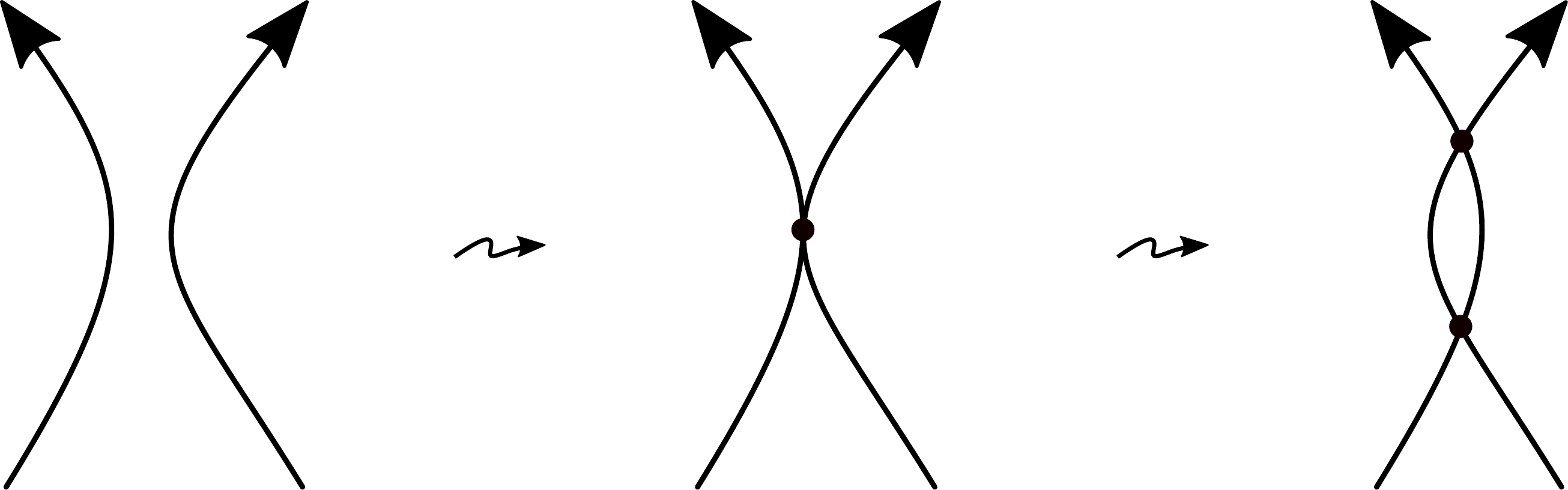}
	\caption{Direct self tangency during a regular homotopy. Left to right: positive. Right to left: negative.}
	\label{fig:jselfdirect}
\end{figure}

\begin{figure}[h!]
	\centering
	\includegraphics[scale=0.3]{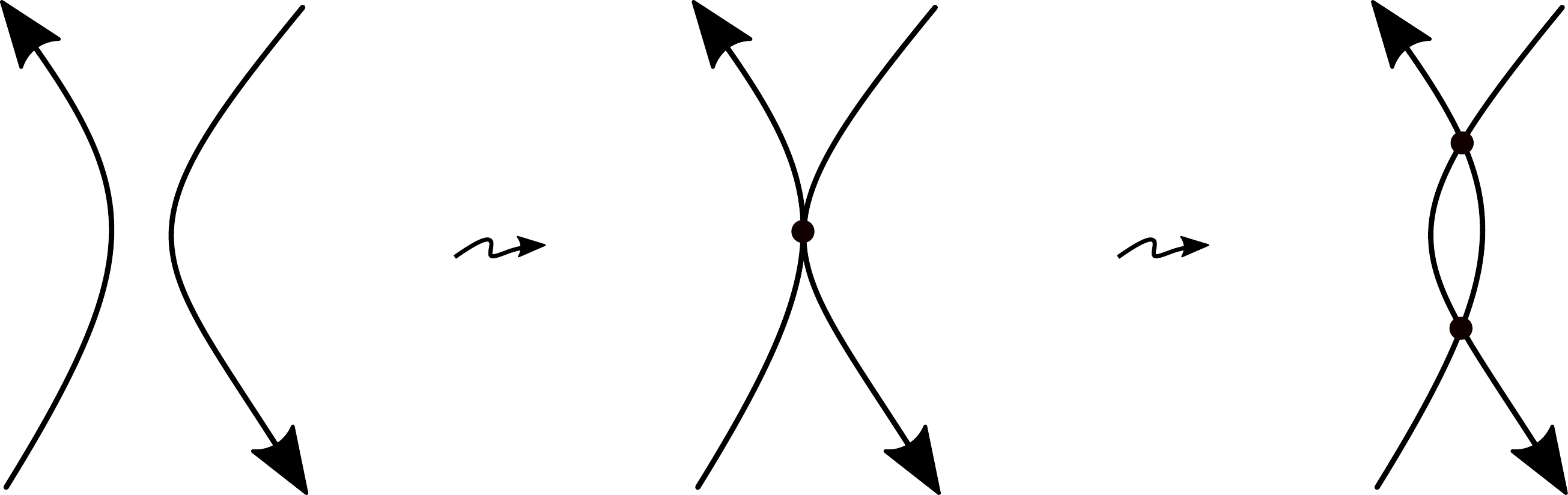}
	\caption{Inverse self tangency during a regular homotopy. Left to right: positive. Right to left: negative.}
	\label{fig:jselfinverse}
\end{figure}

\begin{figure}[h!]
	\centering
	\includegraphics[scale=0.35]{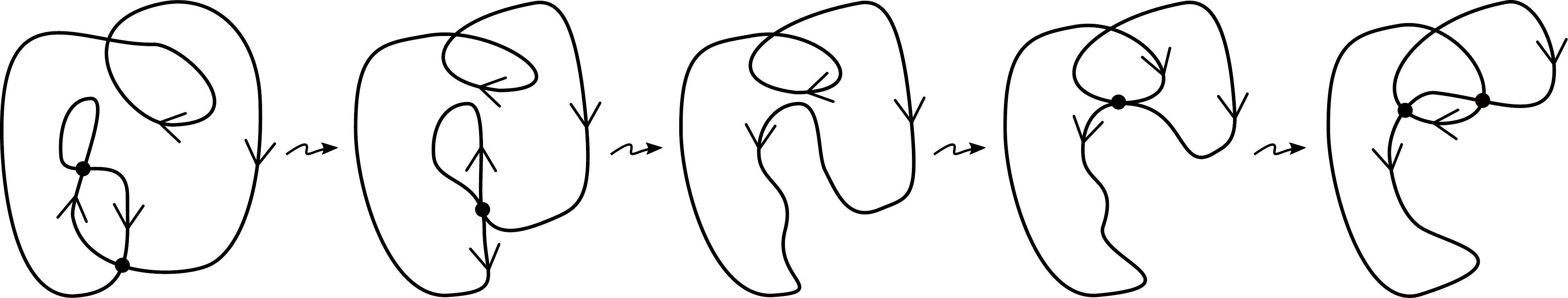}
	\caption{Regular homotopy between two immersed loops with (reading left to right) one negative inverse self-tangency (picture 2) and one positive direct self-tangency (picture 4), first removing, then adding two double points.}
	\label{fig:homtangex}
\end{figure}

\begin{figure}[h!]
	\centering
	\includegraphics[scale=0.3]{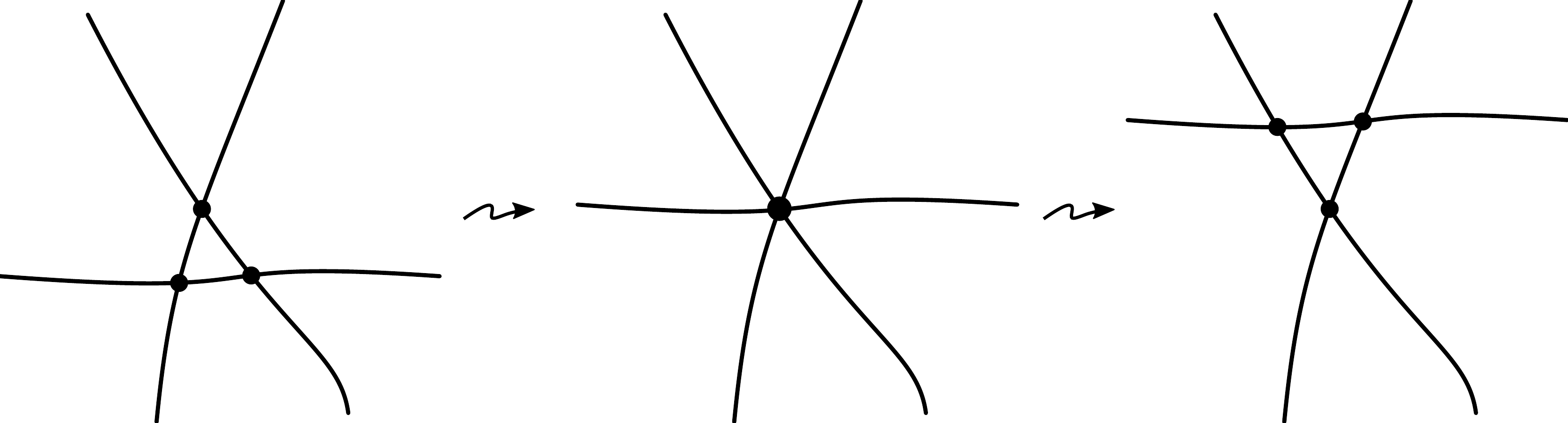}
	\caption{Triple point occurring during a regular homotopy.}
	\label{fig:jtrip}
\end{figure}

Only direct self-tangencies are important for the $J^+$-invariant -- as well as the closely related $\mathcal{J}_1$- and $\mathcal{J}_2$-invariants -- while inverse self-tangencies have no effect on it, but on the $J^-$-invariant. Another event that can occur, that has no effect on any of our invariants, is the crossing of triple points as seen in Figure~\ref{fig:jtrip}. The number of double points does not change at this event.

\subsection{\texorpdfstring{The $J^{\pm}$-invariants}{The J±-invariants}}
\label{subsec:jplusminusinvs}

In his paper \cite[Plane Curves, Their Invariants, Perestroikas and Classifications]{arnold:paper}, Arnold introduces the invariants $J^+, J^-$ and $St$ among others. His findings, especially on the $J^+$-invariant and its well-definedness, lay the foundations for the results of this paper.

\begin{figure}[h!]
	\centering
	\includegraphics[scale=0.4]{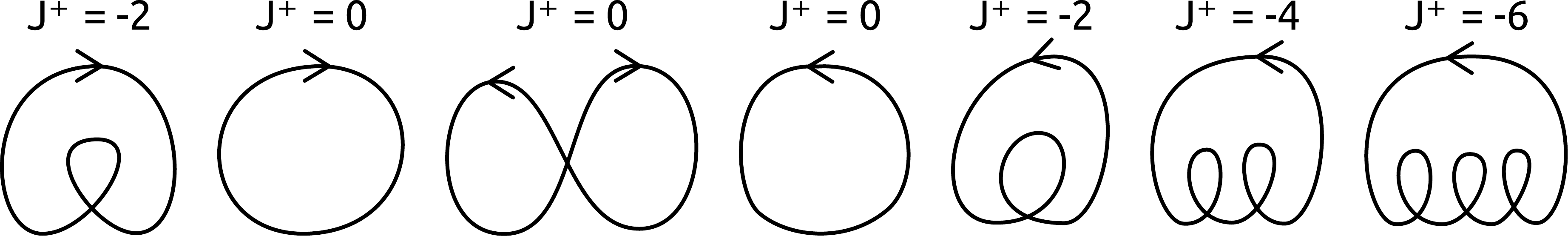}
	\caption{The standard curves. From left to right: $K_{-2}, K_{-1}, K_0, K_1, K_2, K_3, K_4 .$}
	\label{fig:standardcurves}
\end{figure}

\begin{defi}[Standard curves $K_j$]
	\label{def:standardcurves}
	We call the immersions in Figure~\ref{fig:standardcurves} the \emph{standard curves} $K_j .$
	\begin{itemize}
		\item $K_0$ is the figure eight,
		\item $\forall j \neq 0: K_j$ is a circle with $|j| - 1$ many interior single loops that do not intersect with each other and with rotation number~$\operatorname{rot}(K_j) = j .$
	\end{itemize}
\end{defi}

\begin{defi}[$J^+$]
	\label{def:jplus}
	The invariant $J^+$ is a map
	$$\{ K \, | \, \text{$K$ is a generic immersed loop} \} \longrightarrow 2 \mathbb{Z},$$
	with
	$$K_0 \longmapsto 0$$
	and
	$$K_j \longmapsto -2 (|j| - 1) .$$

	Further, the invariant $J^+$
	\begin{itemize}
		\item changes by $+2$ ($-2$) under positive (negative) direct self-tangencies, i.e.\@ tangent immersion crossings where the number of double points increases (decreases) and both involved tangent arcs have the same direction (see Definition~\ref{def:selftang}),
		\item does not change under inverse self-tangencies or crossings through triple points,
		\item is additive under \emph{connected sums},
		\item is independent of the orientation of the immersion.
	\end{itemize}
\end{defi}
\vspace*{0.8em}

\begin{defi}[$J^-$]
	\label{def:jminus}
	The invariant $J^-$ is a map
	$$\{ K \, | \, \text{$K$ is a generic immersed loop} \} \longrightarrow \mathbb{Z},$$
	with
	$$K_0 \longmapsto -1$$
	and
	$$K_j \longmapsto -3 (|j| - 1) .$$

	Further, the invariant $J^-$
	\begin{itemize}
		\item changes by $-2$ ($+2$) under positive (negative) inverse self-tangencies, i.e.\@ tangent immersion crossings where the number of double points increases (decreases) and the two involved tangent arcs have opposite directions (see Definition~\ref{def:selftang}),
		\item does not change under direct self-tangencies or crossings through triple points,
		\item is additive under \emph{connected sums},
		\item is independent of the orientation of the immersion.
	\end{itemize}
\end{defi}
\vspace*{0.8em}

With $n_K$ the number of double points of an immersion $K,$ Arnold shows (see \cite[Page 16]{arnold:paper}), that:
$$J^+(K) - J^-(K) = n_K .$$

So we have:
\begin{equation}
	\label{eq:arnoldjinvsdps}
	J^-(K) = J^+(K) - n_K .
\end{equation}

In Chapter~\ref{subsec:szhom} we will introduce two further invariants $\mathcal{J}_1$ and $\mathcal{J}_2$. Just as $J^+$ and $J^-,$ they are independent of the orientation of the immersion. Most of the time we only denote the orientation in the picture of an immersion to decide whether a self-tangency is direct or inverse.

With this in mind, some pictures of immersions from here on will lack any indication of orientation.

\subsubsection{Viro's formula}%
\label{subsubsec:viro}%

\begin{thesislemma}\label{lem:viro}[Viro's Formula]
	Let $K$ be a generic immersed loop. Then
	$$ J^+(K) = 1 + n_K - \sum \limits_{C \in \Gamma_K} (\omega_C(K))^2 + \sum \limits_{p \in \mathcal{D}_K} (\operatorname{ind}_p(K))^2 ,$$
	with
	\begin{itemize}
		\item $n_K$ the number of double points of $K$
		\item $\Gamma_K$ the connected components of $\mathbb{C} \setminus K$
		\item $\mathcal{D}_K$ the double points of $K$
		\item $\operatorname{ind}_p(K)$ the index of the double point $p$ in $K ,$ see Definition~\ref{def:indexdp}
	\end{itemize}
\end{thesislemma}

\begin{remark}
	In Viro's paper, the formula is stated as
	$$J^+(K) = 1 - \sum \limits_{C \in \Gamma_K} (\omega_C(K))^2 + \sum \limits_{p \in \mathcal{D}_K} (1 + (\operatorname{ind}_p(K))^2) ,$$
	which is equal to the way we wrote it.

	With Equation~\ref{eq:arnoldjinvsdps} we can easily use Viro's formula to calculate $J^-$ by just leaving out the addition of $n_K$:
	$$ J^-(K) = 1 - \sum \limits_{C \in \Gamma_K} (\omega_C(K))^2 + \sum \limits_{p \in \mathcal{D}_K} (\operatorname{ind}_p(K))^2 .$$

	The proof for this formula can be done in several different ways. We will not do any of them. The easiest full proof of Viro's formula, that can be done with only the concepts that we used so far in this paper, can be found on page 7 of \textnormal{The $J^{2+}$-Invariant for Pairs of Generic Immersions} by \textnormal{Hanna Häußler} \cite{hanna:paper}. Of course it is also proven in the cited literature by Oleg Viro \cite{viro:paper}, but differently.
\end{remark}

An example using Viro's formula can be found in the appendix on page~\pageref{example:viro}.

\subsection{\texorpdfstring{$\mathcal{J}_1$- and $\mathcal{J}_2$-invariant and Stark--Zeeman homotopies}{J1- and J2-invariant and Stark--Zeeman homotopies}}
\label{subsec:szhom}

This chapter is not needed for the results about the $J^+$-invariant of $k$-bifurcations of immersions, which we present in Chapter~\ref{subsec:kbifj}. It is only needed for the corollaries in the last two Chapters~\ref{subsec:bifj1} and~\ref{subsec:bifj2}.

In order to classify families of periodic orbits in planar Stark--Zeeman systems, \emph{Cieliebak, Frauenfelder and van Koert} introduced the $\mathcal{J}_1$- and $\mathcal{J}_2$-invariant based on Arnold's $J^+$-invariant, defined for generic closed plane immersions~$K \subset \mathbb{C} \setminus \{ 0 \}$ \cite{kai:paper}.

These two new invariants rely on an origin point as they change depending on the \emph{winding number}~$\omega_0(K) .$ Further, they are invariant under \emph{Stark--Zeeman homotopies} (introduced shortly). While they determine each other for immersions with odd winding number \cite[Proposition 6]{kai:paper} and some other specific types of orbits \cite[Lemma 5]{kai:paper}, the pair $(\mathcal{J}_1(K), \mathcal{J}_2(K))$ attains all values in $2\mathbb{Z} \times 2\mathbb{Z}$ for immersions~$K$ with $\omega_0(K)$ even \cite[Proposition 7]{kai:paper}.

\begin{figure}[h!]
	\centering
	\includegraphics[scale=1.2]{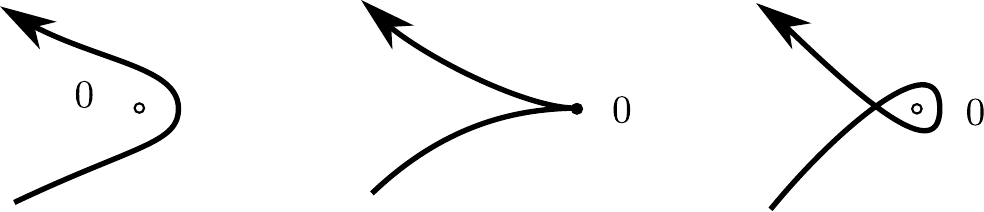}
	\caption{Event $I_0$ (figure taken from \cite[Figure 5]{kai:paper}, slightly adjusted). Left to right: birth of interior loop. Right to left: death of interor loop. Note that in the left picture the curve was going around the origin, i.e.\@ the added loop ended up in the component of the curve that did not include the origin.}
	\label{fig:intloopcusp}
\end{figure}

\begin{defi}[Stark--Zeeman homotopy]
	\label{def:szhom}
	A \emph{Stark--Zeeman homotopy} (from \cite[Definition 1]{kai:paper}) between two generic immersed loops $q$ and $q'$ is a smooth map $$h : S^1 \times [0, 1] \to \mathbb{C},$$ with $h(\cdot, 0) = q$ and $h(\cdot, 1) = q' ,$ and $h(\cdot, t) : S^1 \to \mathbb{C}$ a generic immersed loop $\forall t \in [0, 1]$ except for a finite number of $t \in (0, 1)$ where these isolated events can happen:

	\begin{itemize}
		\item $(I_0)$ birth or death of an interior loop through a cusp at the origin, see Figure~\ref{fig:intloopcusp}
		\item $(I_\infty)$ birth or death of an exterior loop (interior loop of the unbounded component) through a cusp
		\item $(II^+)$ inverse self-tangency (see Figure~\ref{fig:jselfinverse})
		\item $(III)$ triple point crossing (see Figure~\ref{fig:jtrip})
	\end{itemize}
\end{defi}

For the physical motivation of these additional moves see~\cite{kai:paper}.

In other words, a Stark--Zeeman homotopy is a regular homotopy that is not allowed to let the immersion cross the origin and where exterior loops and interior loops around the origin (as seen in Figure~\ref{fig:intloopcusp}) can be added or removed. Note that Stark--Zeeman homotopies change the rotation number of an immersion when the events $I_0$ or $I_\infty$ occur.

Let us denote $\mathbb{C} \setminus \{ 0 \}$ as $\mathbb{C}^*$ from here on.

\begin{defi}[$\mathcal{J}_1$-invariant]
	\label{def:j1}
	The invariant $\mathcal{J}_1$ is the map
	$$ \{ K \, | \, \text{$K$ is a generic immersed loop in $\mathbb{C}^*$} \} \longrightarrow 2 \mathbb{Z} \cup (2 \mathbb{Z} + \textstyle\frac{1}{2}) ,$$
	$$ K \longmapsto J^+(K) + \dfrac{\omega_0(K)^2}{2} .$$
\end{defi}

\begin{defi}[$\mathcal{J}_2$-invariant]
	\label{def:j2}
	The invariant $\mathcal{J}_2$ is the map
	$$ \{ K \, | \, \text{$K$ is a generic immersed loop in $\mathbb{C}^*$} \} \longrightarrow 2 \mathbb{Z} ,$$
	$$ K \longmapsto \begin{cases}
		J^+(L^{-1} (K)), &\text{if \( \omega_0(K) \) odd} \\
		J^+(\widehat{K}), &\text{else}
	\end{cases} $$
	with $L: \mathbb{C}^* \to \mathbb{C}^*, v \mapsto v^2$ the complex squaring map and $\widehat{K}$ one of the two up to rotation identical immersions of the preimage $L^{-1} (K) .$
\end{defi}
\vspace*{0.8em}

The notion of \emph{lifts} is recalled in the next subchapter, which also elaborates on the preimage of the squaring map $L .$ See the appendix at the end of this paper (page~\pageref{website:fibers}) for an interactive website that visualizes the preimage of $L .$

For the case $\omega_0(K)$ odd, the preimage $L^{-1} (K)$ consists of two point-symmetric (with respect to the origin) curves that are not closed. The union of these two curves is a generic closed immersion, which we denote by $L^{-1} (K) .$

For the case $\omega_0(K)$ even, the preimage $L^{-1} (K)$ consists of two point-symmetric (with respect to the origin) generic closed immersions. The union of these two immersions is not an immersion. Because they are identical up to rotation, their $J^+$-value is the same, so we pick any of the two and denote it by~$\widehat{K} .$

In \cite{kai:paper} it is written as $\widetilde{K}$ instead of $\widehat{K} ,$ which we change in this paper, as we write $\widetilde{K}_k$ for $k$-bifurcations of the immersion~$K$ and want to avoid confusion. The hat above $K$ pointing upwards is chosen to reflect its role as the lift of $K ,$ as lifts are often associated as ``lying above'' in visualizations of lifts.

In their paper, \emph{Cieliebak, Frauenfelder and van Koert} show that $\mathcal{J}_1$ and $\mathcal{J}_2$ are both invariant under Stark--Zeeman homotopies \cite[Propositions 4 \& 5]{kai:paper}. Note that they are also invariant under regular homotopies as long as the immersion does not cross the origin during the homotopy.

Later we will calculate $\mathcal{J}_2$ for bifurcations of immersions. To do this, we need to take a look at the double points of the preimage $L^{-1} (K) .$

\subsubsection{\texorpdfstring{Lifts of the squaring map~$L$ and their double points}{Lifts of the squaring map L and their double points}}
\label{subsubsec:liftsdpsj2}

Please note: the following observations are only of interest for the results on $\mathcal{J}_2$ and can otherwise be ignored.

\begin{figure}[h!]
	\centering
	\includegraphics[scale=0.4]{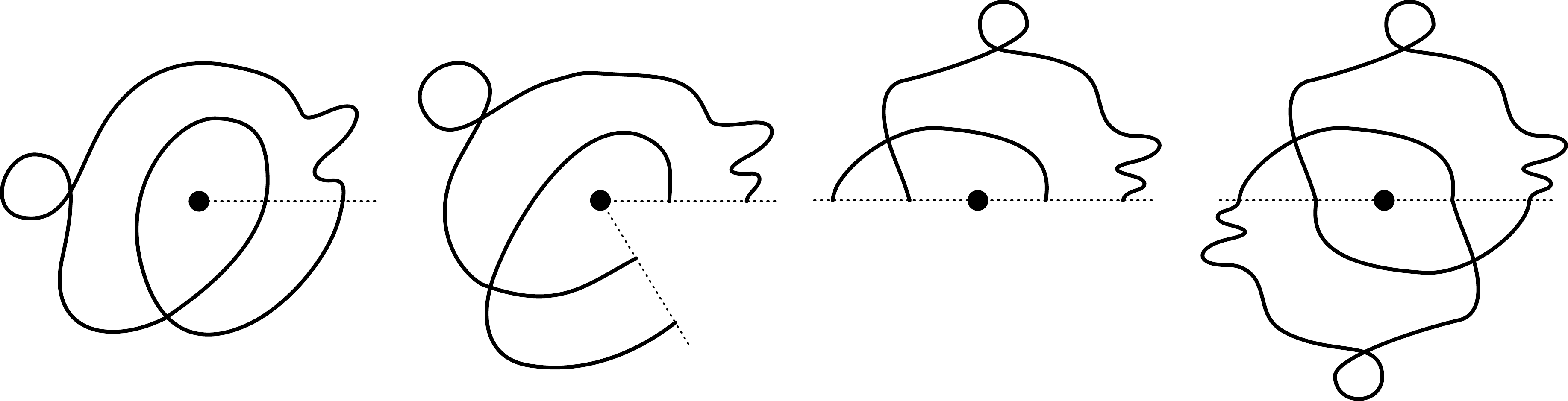}
	\caption{First picture is the immersion with the cut on the positive real axis. Second picture is the ``folding fan'' closed by a little, third picture half closed. In the last picture it is duplicated and one of them rotated by $180^\circ$ around the origin. Here we ignored the squaring of the absolute value of all points. The appendix at the end of this paper links to a website where preimages of the squaring map can be easily experimented with. Note that the figures in this subchapter are hand-drawn, not computer-generated, so not 100\% accurate.}
	\label{fig:squaringinverserip}
\end{figure}

The complex squaring map~$L$ squares the absolute value of a point in $\mathbb{C}$ and doubles its angle in polar coordinates. So to imagine what the preimage of $L$ looks like, ignore the squaring of the absolute value -- this is not dishonest, as this does not change the topological properties of an immersion -- and interpret the reverse of the doubling of the angle like this:

Cut the complex plane $\mathbb{C}$ from the origin~$0$ along the positive real axis. Now imagine this cut open plane as a fully opened folding fan and evenly close it halfway. This halves all the angles. As $L$ is a two-to-one mapping, we need to duplicate this halfway closed folding fan, turn it $180^\circ$ around the origin and stick the two fans together. This process is illustrated in Figure~\ref{fig:squaringinverserip}.

Let us recall some language from (algebraic) topology first and then study which double points of $K$ appear in its preimage. We will need this to classify which double points of an immersion~$K$ we need to consider for the calculation of $\mathcal{J}_2(K)$ in Chapter~\ref{subsec:bifj2}.

We call the squaring map~$L$ a \emph{projection} from the \emph{covering space}~$\mathbb{C}^*$ onto the \emph{base space}~$\mathbb{C}^* .$ We know that $L$ is continuous, even smooth and locally diffeomorphic, because squaring a complex number is the same as doubling its angle in polar coordinates and squaring its absolute value, which are both smooth and locally diffeomorphic maps. Further, $L$ is a \emph{double cover} of $\mathbb{C}^*$ since $L$ is a two-to-one mapping of $\mathbb{C}^*$ onto itself.

\begin{figure}[h!]
	\centering
	\includegraphics[scale=0.5]{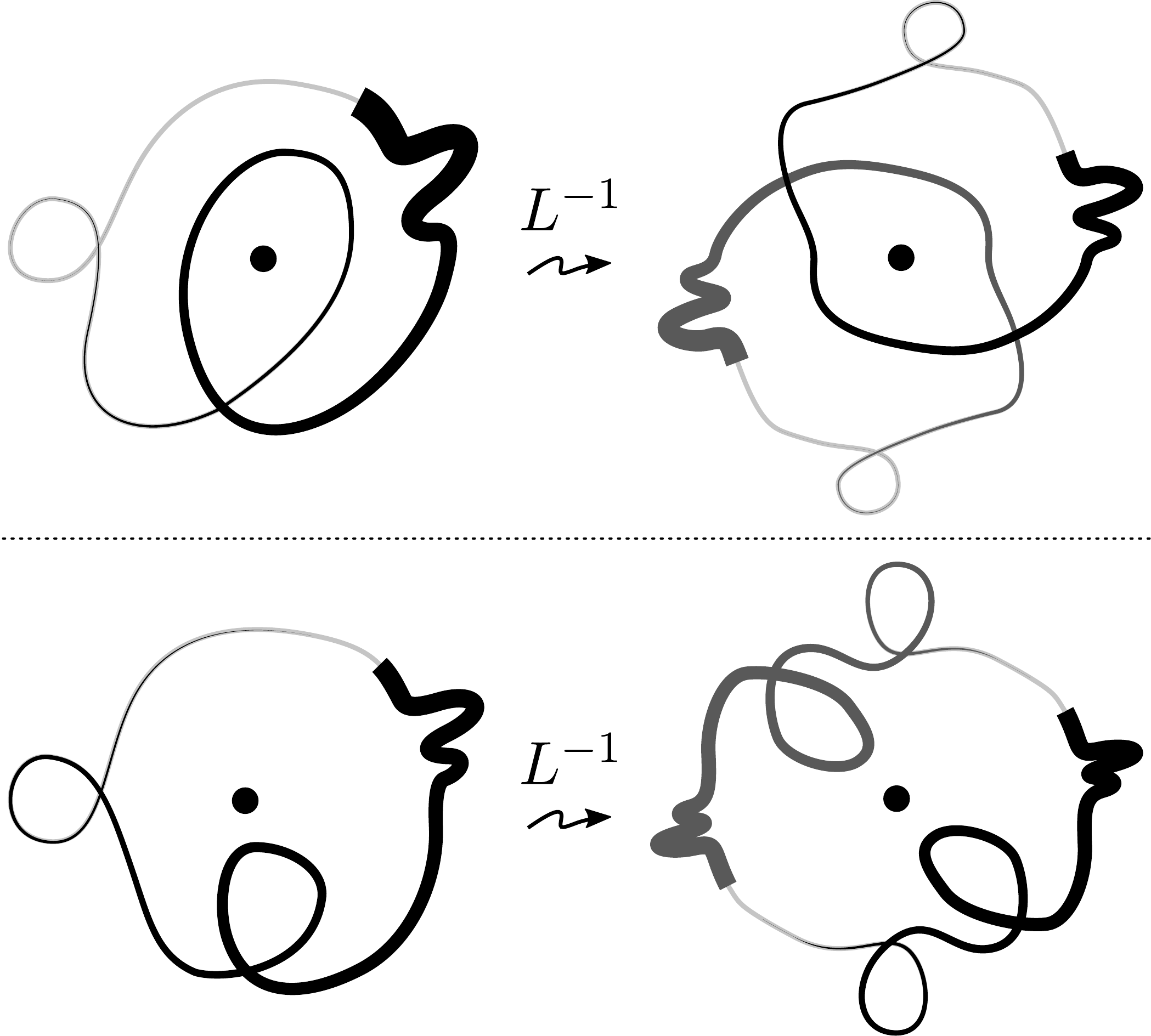}
	\caption{On the left an immersion and on the right its preimage for the squaring map. Top row shows an immersion with even winding number, bottom row with odd winding number. Note that the figures in this subchapter are hand-drawn, not computer-generated, so not 100\% accurate.}
	\label{fig:preimagetwocurves}
\end{figure}

When we look at the preimage of the points that we would visit in one run along the immersion~$K ,$ we see that the preimage draws two immersions simultaneously, see Figure~\ref{fig:preimagetwocurves} for two examples.

If the winding number~$\omega_0(K)$ is odd, these two immersions in the preimage by themselves are not closed immersions. But together they form a single closed immersion, which we simply denoted as~$L^{-1} (K)$ above. Note that the two immersions in the preimage differ depending on where we start our run along $K ,$ but their union~$L^{-1} (K)$ is always the same.

If the winding number~$\omega_0(K)$ is even, these two immersions in the preimage by themselves are closed immersions that are identical up to rotation by $180^\circ .$ Here it does not matter at which point we started our run along $K ,$ as the two immersions in the preimage are closed.

For both cases, we call each of the two immersions in the preimage -- which are not closed in the case~$\omega_0(K)$ odd -- a \emph{lift} of $K .$ There are no other covering maps in this paper, so assume $L$ to be the covering map whenever we talk about lifts of an immersion. This concept is only used again in Chapter~\ref{subsec:bifj2}.

We also take note that due to the properties of the squaring map~$L$ the lifts of an immersion~$K$ are always identical up to rotation by $180^\circ .$ From this we know that any double point in one of the lifts also appears in the other lift. \\

Next let us study which double points of an immersion~$K$ appear in its preimage. Let $n_K$ be the number of double points of $K .$ Because the squaring map~$L$ is a locally diffeomorphic two-to-one mapping, it is easy to see that all double points of $K$ appear twice in the preimage of $K .$ For the case that $\omega_0(K)$ is odd this means that the preimage~$L^{-1} (K)$ has $2n_K$ many double points.

But remember that $\mathcal{J}_2(K)$ of any immersion $K$ with even winding number is just $J^+(\widehat{K}) ,$ with $\widehat{K}$ one of the two up to rotation identical lifts of $K .$ So we want to know which double points of $K$ are double points in the lifts of $K ,$ as intersections between the two lifts are ignored for the calculation of~$J^+(\widehat{K}) .$ We need to distinguish between two different types of double points in $K$ with even winding number~$\omega_0(K) .$

Consider an arbitrary double point~$p$ of an arbitrary immersion~$K$ with winding number~$\omega_0(K)$ even. As always in this paper, the double points in $K$ are all transversal, so no tangential self-intersections. Now let us decompose the immersion~$K$ into two closed curves~$K_p^A$ and~$K_p^B$ as follows.

\begin{figure}[h!]
	\centering
	\includegraphics[scale=0.33]{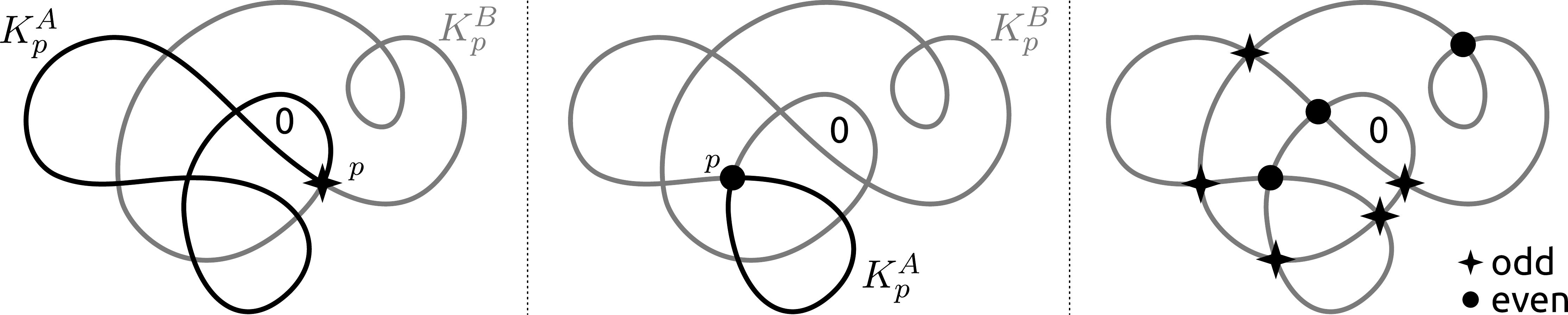}
	\caption{The first two pictures show an immersion with even winding number -- note the $0$ in the picture indicating the origin -- and their decompositions at a point $p .$ The third picture is the same immersion, but with all double points labeled whether they are even or odd, according to Definition~\ref{def:evendp}.}
	\label{fig:imdecompatp}
\end{figure}

Consider the four immersion arcs that meet at the double point~$p$ -- two are directed away from~$p ,$ two are not. We define $K_p^A$ as the closed curve that we get if we start at $p$ and follow one of the two immersion arcs that are directed away from $p$ until we return to $p$ for the first time. Note that $K_p^A$ is generic, but not regular or even smooth at the point $p ,$ which is also not a double point in $K_p^A .$ We define $K_p^B$ as the other closed curve that we can get by choosing the other immersion arc that is directed away from $p .$ Figure~\ref{fig:imdecompatp} shows examples of such decompositions.

Remember that the immersion~$K$ has even winding number~$\omega_0(K) ,$ so the winding numbers~$\omega_0(K_p^A)$ and~$\omega_0(K_p^B)$ around the origin of the closed curves~$K_p^A$ and~$K_p^B$ are either both even or both odd. With this we define the following -- where of course $K_p^A$ can be replaced with $K_p^B .$

\begin{defi}[Even (odd) double point]
	\label{def:evendp}
	Let $p$ be an arbitrary double point of an arbitrary immersion~$K$ with even winding number~$\omega_0(K)$ around the origin~$0$ and $K_p^A$ as described above and illustrated in Figure~\ref{fig:imdecompatp}.

	We define $p$ to be an \emph{even (odd) double point} if the winding number~$\omega_0(K_p^A)$ is even (odd).
\end{defi}

Now we can easily see that only even double points of $K$ appear as double points in a lift~$\widehat{K} .$ If $p$ is an odd double point, then a lift of $K_p^A$ is not a closed curve and neither its start point nor its end point is a double point in $\widehat{K} .$ Instead, the preimage of an odd double point is the set of two double points where the two lifts of $K$ intersect. This can also be seen in the top row of Figure~\ref{fig:preimagetwocurves}.

\begin{thesiscorollary}
	\label{cor:evenliftdps}
	Let $K$ be an arbitrary immersion with winding number $\omega_0(K)$ even, $n_K$ its number of double points and $\nu_K$ its number of even double points. Let $\widehat{K}^1$ and $\widehat{K}^2$ be the two up to rotation identical lifts of $K .$ Further let $n_{\widehat{K}^1 \times \widehat{K}^2}$ be the number of intersections between the two lifts. Then:
	$$\nu_K = n_{\widehat{K}^1} = n_{\widehat{K}^2}, \quad n_K - \nu_K = \textstyle\frac{1}{2} n_{\widehat{K}^1 \times \widehat{K}^2}$$
\end{thesiscorollary}

Later in Chapter~\ref{subsec:bifj2} we will use these observations to calculate the change of $\mathcal{J}_2$ under $k$-bifurcations of immersions.

\clearpage

\section{Invariants under bifurcations}
\label{sec:jbif}

For any immersion~$K$ -- as before we only talk about generic immersed loops and simply call them immersions -- a lower bound for $J^+$ is given by \emph{Arnold's Conjecture}~\cite[Page 33]{arnold:paper}, which was later proved by Viro~\cite[Conjecture 1.3.A]{viro:paper}:
\begin{equation}
	\label{eq:arnoldsconjecture}
	J^+(K) \ge -n^2 - n ,
\end{equation}
with $n$ the number of double points of $K.$

In this chapter, we will first discuss different bifurcations of immersions and then, using Viro's formula (see Lemma~\ref{lem:viro}), we give a far better lower bound for bifurcations, showing that for any $k$-bifurcation~$\widetilde{K}_k$ of any generic immersed loop $K$ the following attainable inequality holds:
$$ J^+(\widetilde{K}_k) \ge k^2 J^+(K) - (k^2 - k) .$$

Of course for all appearances of $k$ in this chapter we have $k \in \mathbb{N}, \, k \ge 2 .$

In addition, the inequality is not only a lower bound, but attained for each immersion's $k$-bifurcation with no additional crossings, a notion to be introduced below, equivalent to minimal number of double points. These $k$-bifurcations with zero additional crossings are constructed and proven to be of minimal~$J^+.$

From this we deduce a wealth of corollaries. Not only for $J^+$ but also for $J^-$, $\mathcal{J}_1$ and $\mathcal{J}_2 .$

But first we discuss bifurcations and why $J^+(\widetilde{K}_k)$ has no upper bound -- as far as bifurcations are understood in this paper.

\subsection{\texorpdfstring{Visualizing $k$-bifurcations}{Visualizing k-bifurcations}}
\label{subsec:visbif}

To get a $k$-bifurcation of an immersion $K \subset \mathbb{C}$ we run through the immersion $k$ times and then perturb it to a generic immersion, so all self-intersections are transverse double points. We assume the perturbations to be very small and that they do not result in intersections between runs of two parts of the original immersion that do not intersect.

\begin{defi}[Path and strands]
	\label{def:pathstrands}
	When talking about $k$-bifurcations of an immersion~$K$ we denote a sufficiently small neighborhood of $K$ as the \emph{path} -- so a thin tube around the pre-bifurcation immersion.
	And denote the $k$ many arcs in the path as \emph{strands}.
\end{defi}

\begin{figure}[h!]
	\centering
	\includegraphics[scale=0.35]{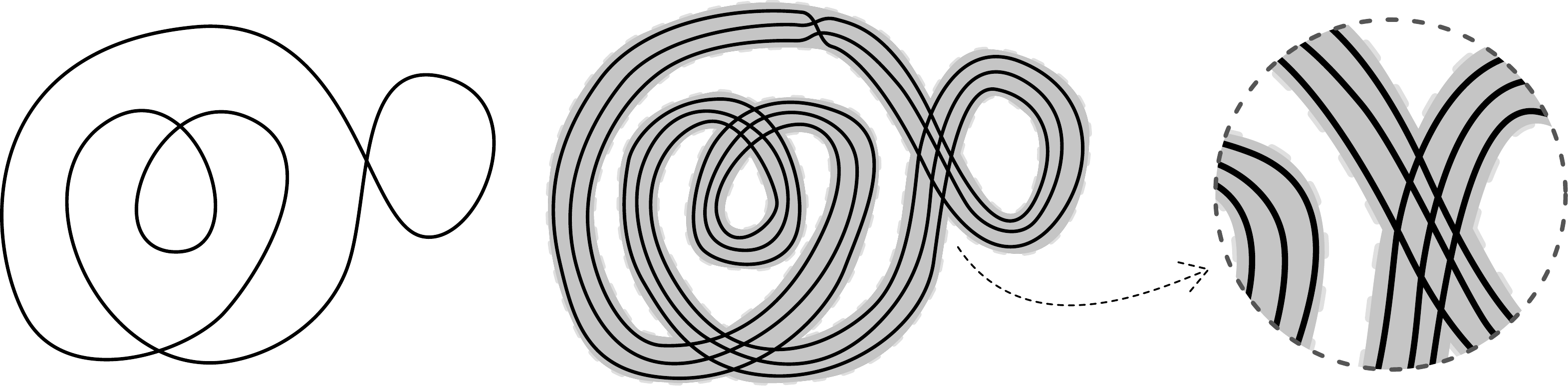}
	\caption{Example of a $3$-bifurcation with close-up of the strands and path for better visibility.}
	\label{fig:pathstrands}
\end{figure}

\begin{remark}
	Note that the concept of strands only makes sense in small neighborhoods of the path that do not contain the whole immersion, as they are all part of a single connected immersion. These denotions of paths and strands are only necessary to point at the right connected components and double points later, not for any rigorous calculations.
\end{remark}

In this paper, when we talk about an immersion [that experienced / after / under] a $k$-bifurcation, we actually mean the immersion at some moment~$\mu \in (\mu_0, \mu_0 + \varepsilon),$ the time interval just after the bifurcation occured and before the immersion experiences further changes.

\begin{figure}[h!]
	\centering
	\includegraphics[scale=0.3]{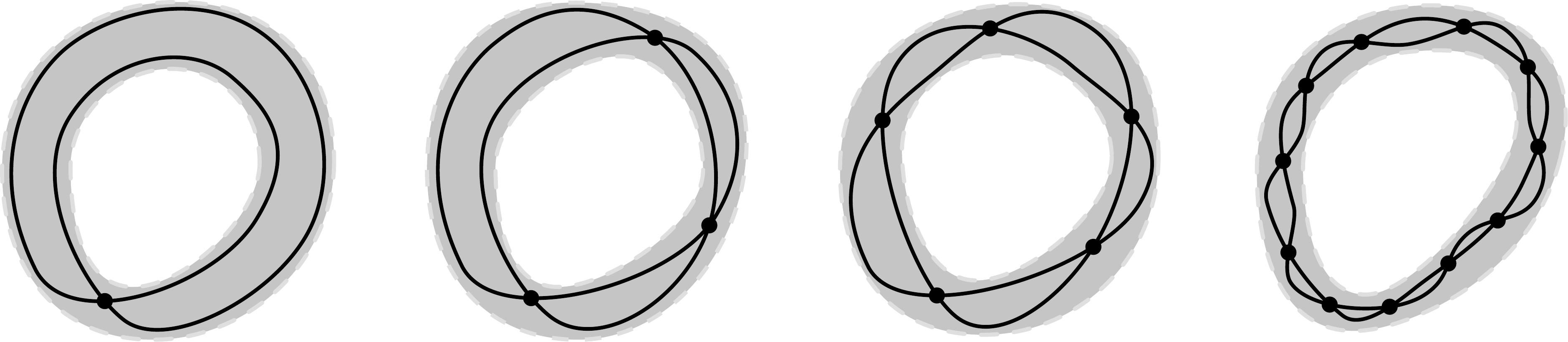}
	\caption{A few $2$-bifurcations of a circle, with visualized paths that unsurprisingly look like circles.}
	\label{fig:moebprojectionspaths}
\end{figure}

\begin{figure}[h!]
	\centering
	\includegraphics[scale=0.3]{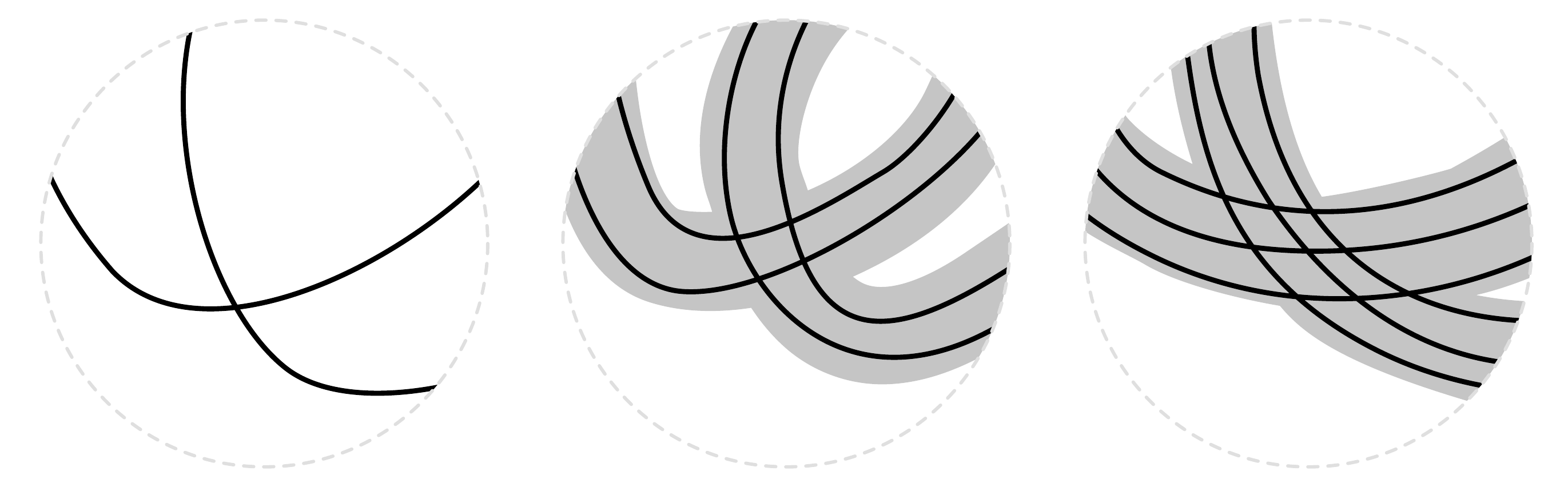}
	\caption{Intersecting path.}
	\label{fig:pathintersect}
\end{figure}

Consequently, the visualizations of bifurcated immersions are always over-dramatizations, as the~$k$ many strands in the path are in general too close together to visualize faithfully when drawing the whole immersion. This is not dishonest, as long as $J+$ and the winding number~$\omega_x(K)$ for an observed point $x \in \mathbb{C} \setminus K$ are not changed. The value of $J^+$ stays the same as long as we ensure that strands within one part of the path do not intersect strands of another part of the path when drawing them, as seen in Figure~\ref{fig:moebprojectionspaths}. Except for when the path intersects itself, as seen in Figure~\ref{fig:pathintersect}.

\begin{figure}[h!]
	\centering
	\includegraphics[scale=0.3]{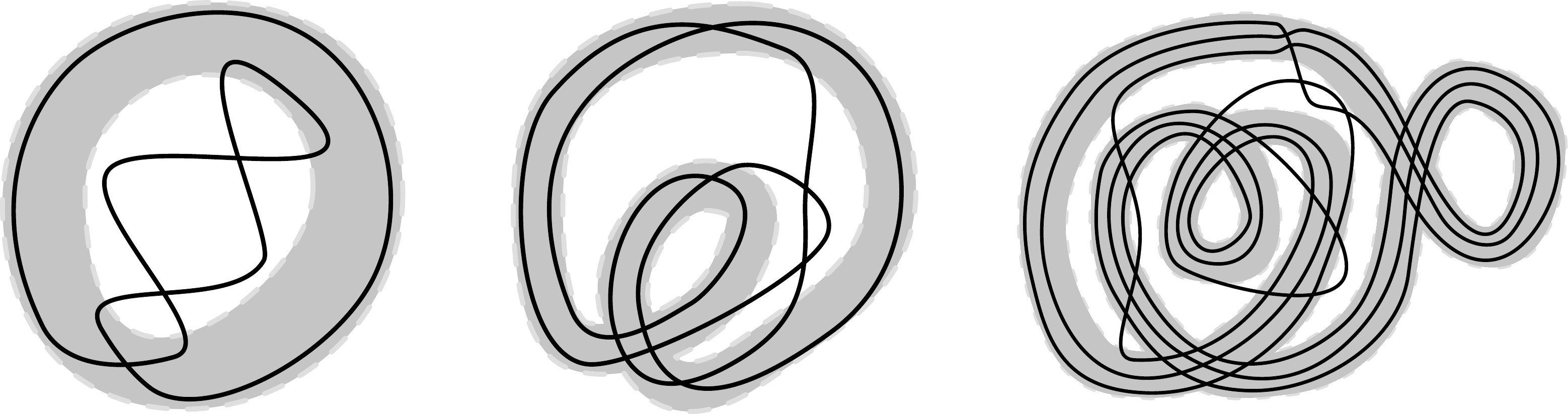}
	\caption{This should not be done.}
	\label{fig:badpathstrands}
\end{figure}

Figure~\ref{fig:badpathstrands} illustrates how strands of a bifurcation should definitely \emph{not} be visualized, even if the pre-bifurcation immersion passed itself very closely, as it shows intersections that should not occur.

Of course if the pre-bifurcation immersion had intersections, the path intersects itself as well and strands of one part of the path intersect all the strands of the other part of the path.

\begin{figure}[h!]
	\centering
	\includegraphics[scale=0.3]{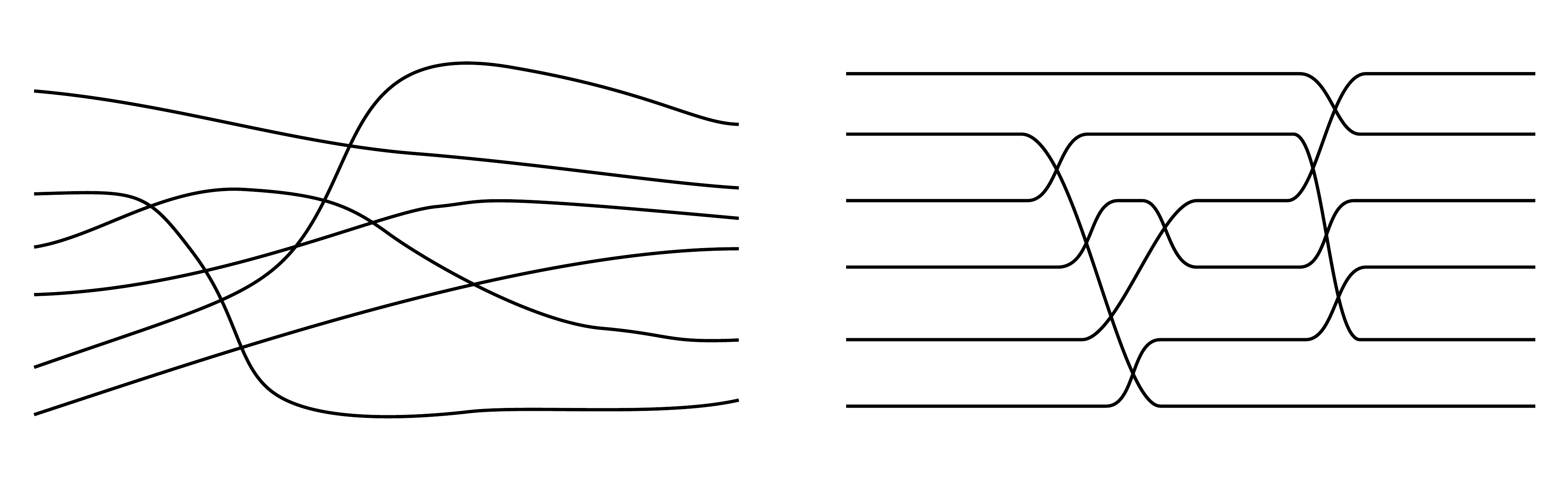}
	\caption{A neighborhood of strands of a bifurcation. On the left it is drawn carelessly, on the right in parallel.}
	\label{fig:parallelstrands}
\end{figure}

We draw strands of a bifurcation roughly parallel to each other whenever possible, especially later when looking at double points between them, to keep pictures tidy and easy to discuss. See Figure~\ref{fig:parallelstrands}.

\subsection{\texorpdfstring{Differences in $k$-bifurcations}{Differences in k-bifurcations}}
\label{subsec:differentkbif}

Let us take for example a circle, i.e.\@ the immersion~$id: S^1 \to S^1 \subset \mathbb{C}, x \mapsto x, $ under a $2$-bifurcation.

\begin{figure}[h!]
	\centering
	\includegraphics[scale=0.3]{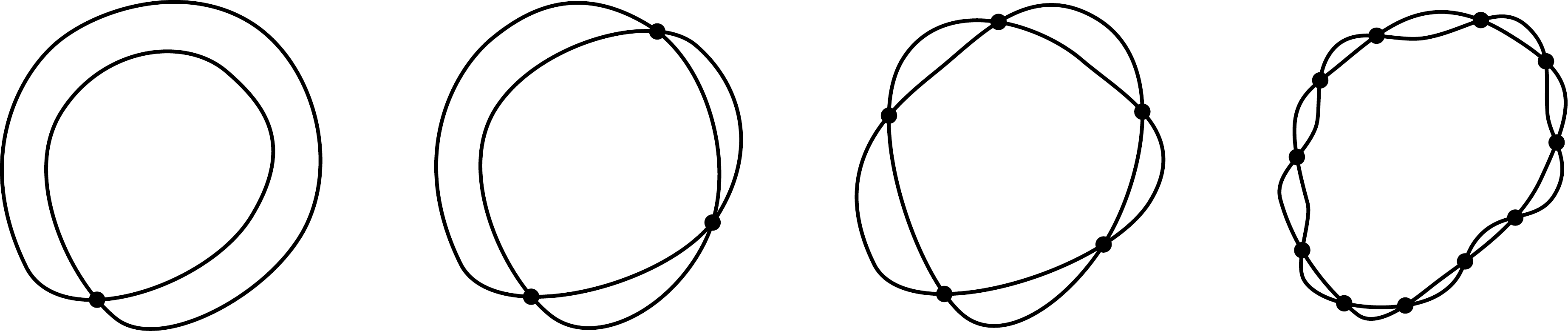}
	\caption{Different $2$-bifurcations of the same circle-shaped immersion.}
	\label{fig:moebprojections}
\end{figure}

When we run through the circle twice and then perturb it to be a general immersion, we can get an arbitrary amount of ``wiggling'' around the original path of the circle. This results in different numbers of double points in the bifurcation, as seen in Figure~\ref{fig:moebprojections}.

\subsubsection{Additional crossings}
\label{subsubsec:addcross}

It is easy to see that there could be an arbitrary even number of added double points in the bifurcation of an immersion.

This is critical to the $J^+$-invariant, as all strands in a part of the path have the same orientation and therefore every single one of these additional crossings -- each creating two new double points through a direct self-tangency, compared to an immersion without it -- increases $J^+$ by $2 .$ Any immersion that experiences a $k$-bifurcation is susceptible to these additional crossings. In Figure~\ref{fig:moebprojections} it is obvious that the first picture is the one with the least additional crossings of the pictured immersions and that it has the lowest $J^+.$

\begin{figure}[h!]
	\centering
	\includegraphics[scale=0.3]{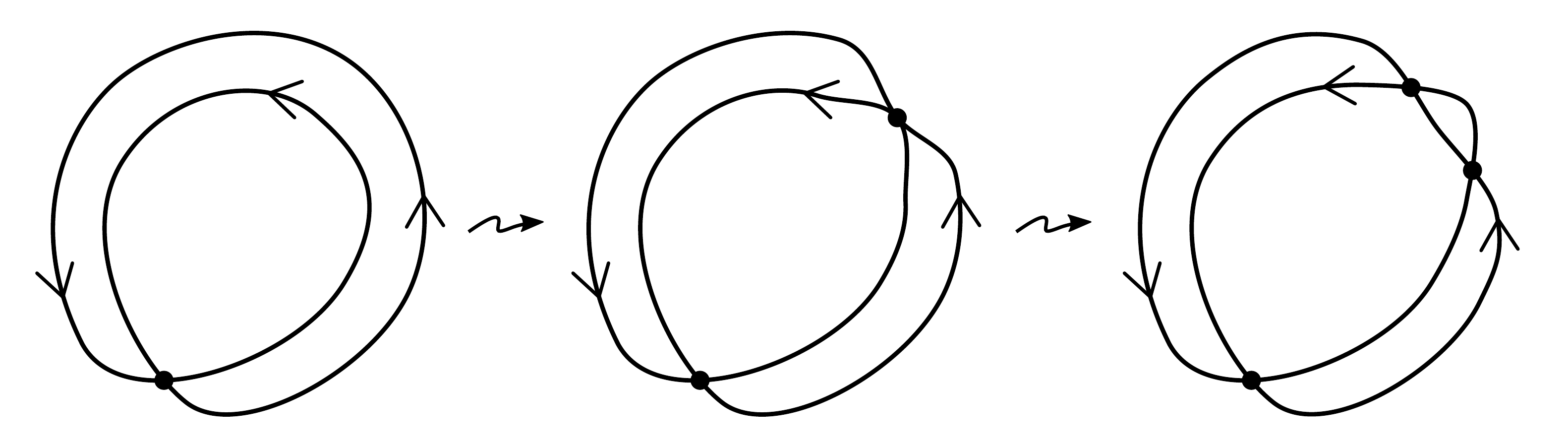}
	\caption{Regular homotopy between two $2$-bifurcations of the circle with different $J^+$ as seen in Figure~\ref{fig:moebprojections}.}
	\label{fig:moebaddcrossing}
\end{figure}

We will call each such crossing between two neighboring strands of an immersion's bifurcation an \emph{additional crossing}. This is to remind us that there is always a regular homotopy that can remove two double points through a negative self-tangency between two neighboring strands, which would change $J^+$ by $-2.$

\begin{remark}
	This uncertainty of how many crossings a $k$-bifurcation has, whether some strands are wavy in a way that creates additional crossings, is the reason why we denote $k$-bifurcations of $K$ with a ``wave'' above~$K ,$ so $\widetilde{K}_k .$
\end{remark}

\begin{remark}
	As a quick side note for later, remember that $J^-$ does not change by additional crossings, as they are always direct self-tangencies, not inverse self-tangencies.
\end{remark}

\subsubsection{\texorpdfstring{Bifurcations with minimal $J^+$}{Bifurcations with minimal J+}}
\label{subsubsec:bifminj}

\begin{figure}[h!]
	\centering
	\includegraphics[scale=0.3]{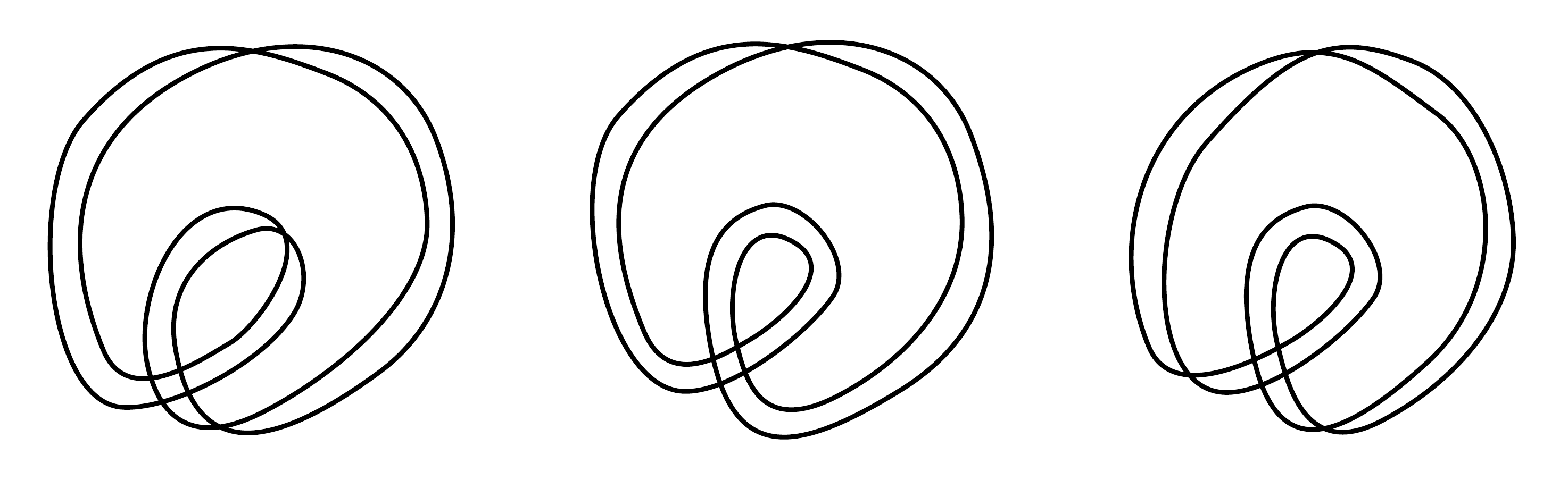}
	\caption{$2$-bifurcations of the standard curve $K_2$.}
	\label{fig:unclearminimalj}
\end{figure}

Without any further observations, it is not immediately clear on first sight which of the immersions in Figure~\ref{fig:unclearminimalj} have the least $J^+$ or even the lowest number of double points.

Maybe intuitively and unsurprisingly, bifurcations with fewer double points should generally have a lower $J^+$-value. Of course without breaking our rules on how to draw bifurcations, see Chapter~\ref{subsec:visbif}. We now show that this intuition is indeed true for double points between strands within the path.

\begin{remark}
	In general the number of double points can be useful for a weak lower bound, as in Arnold's conjecture (see Equation~\ref{eq:arnoldsconjecture}), and is an important part of Viro's formula (see Lemma~\ref{lem:viro}). Figure~\ref{fig:3dpimms} reminds us that the number of double points alone is not enough information to deduce the value of $J^+.$

	\begin{figure}[h!]
		\centering
		\includegraphics[scale=0.3]{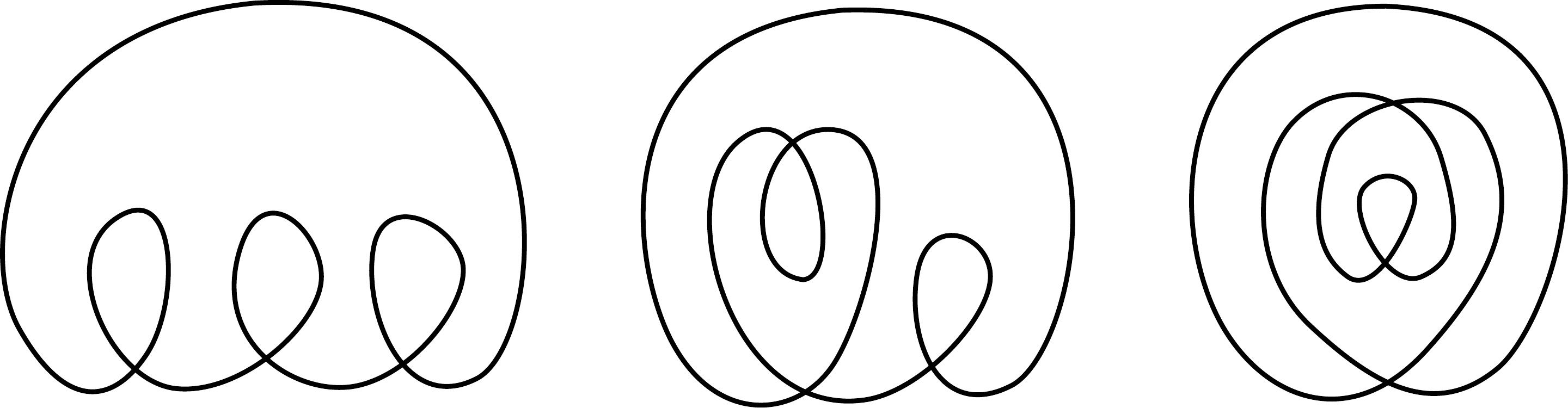}
		\caption{Immersions with three double points each and different $J^+$-value, left to right: $-6, -8, -12.$}
		\label{fig:3dpimms}
	\end{figure}
\end{remark}

Let us first look at the minimal number of double points of a bifurcation.

\begin{thesislemma}
	\label{lem:mindps}
	If $K$ is an immersion with $n_K$ double points, then any $k$-bifurcation $\widetilde{K}_k$ has at least
	$$ n_K k^2 + (k-1) $$
	double points.
\end{thesislemma}

\begin{figure}[h!]
	\centering
	\includegraphics[scale=0.3]{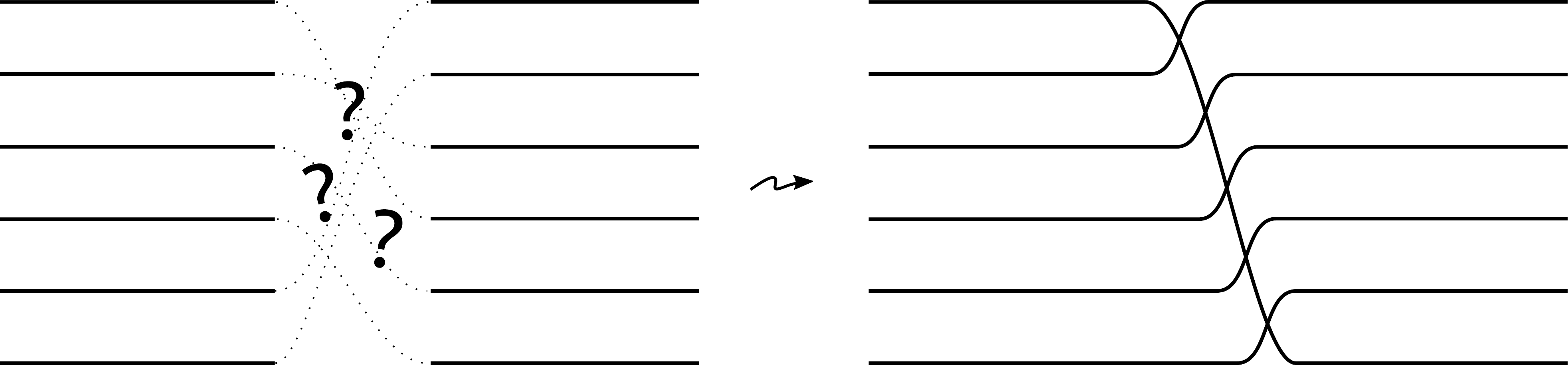}
	\caption{All strands are connected neatly. We now assume that there are no double points between strands at any part of the path anymore, except for exactly $k^2$ double points per path intersection. Here $k = 6 .$}
	\label{fig:neatconnections}
\end{figure}

\begin{proof}
	\label{pr:mindps}
	Let us explain the two parts of the sum separately.

	There are $k^2$ double points in $\widetilde{K}_k$ for every single double point in $K.$

	A double point in $K$ is created by the immersion intersecting itself transversally -- remember that generic immersions have no tangential double points. Without loss of generality, assume that for a small neighborhood of the original double point the $k$ many strands in $\widetilde{K}_k$ run roughly in parallel, i.e.\@ without intersecting each other.

	As Figure~\ref{fig:doublepointbif} illustrates, in the bifurcation $k$ many strands from one part of the immersion each cross~$k$ strands from another part of the immersion, forming a square grid of $k^2$ double points.

	\begin{figure}[h!]
		\centering
		\includegraphics[scale=0.32]{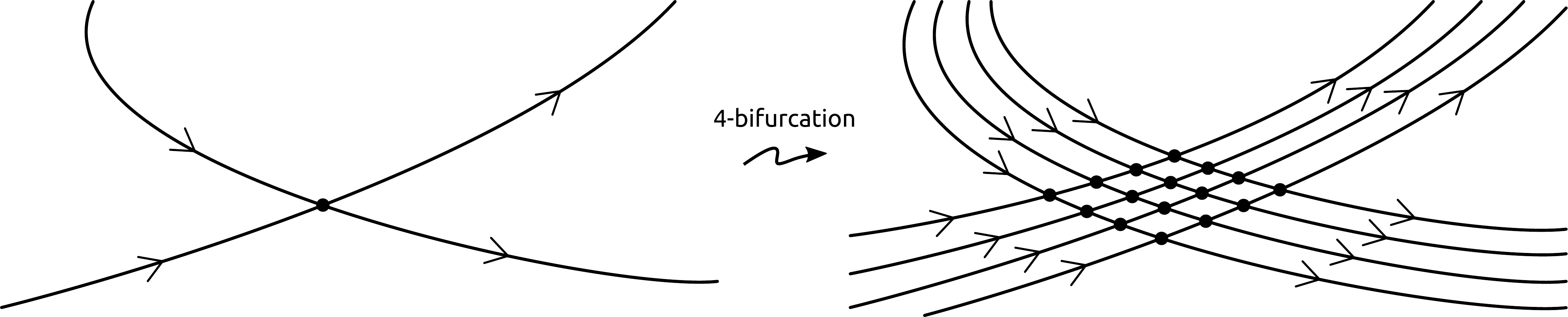}
		\caption{Double points of a $4$-bifurcation in place of a former single double point pre-bifurcation.}
		\label{fig:doublepointbif}
	\end{figure}

	Next, there are at least $k-1$ double points from connecting the $k$ many strands of the bifurcation.

	Let us assume that we draw the $k$ strands of the bifurcation running in parallel the whole way and leave a small part empty where we need to connect the strands together so they form one single bifurcated immersion. Obviously, so far we only have the $n_K k^2$ double points from above.

	To connect the strands all together, we can number them from top to bottom, $1$ to $k ,$ and just take any permutation $\pi \in \mathcal{S}_k ,$ with $\mathcal{S}_k$ the symmetric group, that has \emph{exactly one cycle} so all of them are connected in one immersion. Note that we count elements with $\pi (i) = i$ for any $i \in \{1, \dots, k \}$ as cycles with one element, so if a permutation has one cycle then all elements are in that cycle. The simplest permutation with one cycle would be
	$$ \left(\begin{smallmatrix}
		1 & 2 & 3 & \cdots & k-1 & k \\
		k & 1 & 2 & \cdots & k-2 & k-1
	\end{smallmatrix} \right) ,$$
	which is pictured in Figure~\ref{fig:neatconnections}. We denote this specific connection of strands as \emph{neat strand connections}. It can be written as permutation in cycle notation as $(k, 1, 2, \dots, k - 1) .$
\end{proof}

\begin{remark}[On connections with minimal double points.]
	Note that there are permutations with one cycle that would have more double points in the strand connections and that there are a total of $2^{(k-2)}$ unique permutations with minimal double points -- a proof sketch follows.

	The number of double points for permutations with one cycle is $\left( \sum \limits_{i=0}^{k}  |\pi(i) - i| \right) - (k-1)$ and the minimal amount $\min \limits_\pi \{ \left( \sum \limits_{i=0}^{k}  |\pi(i) - i| \right) - (k-1) \},$ with $\pi$ has one cycle, is $k - 1 .$ The number of double points in the strands is equal to the concept of the \emph{inversion number} for permutations.

	The equation for the number of double points and the result for the minimal amount comes from this simple observation: to connect all strands -- i.e.\@ to get one cycle -- we can start from the bottom strand and see that the connections have to make it all the way to the top-most strand or else it would not be in the cycle. And we have to get all the way back down. If we take no detour on our way up -- this means only going up, irrespective of the connection lengths -- then there are no new double points. But on the way back down $k-1$ double points emerge if we take no detour -- this means only going down and necessarily leaving out no strand that is not yet connected. From this observation we can easily see that the minimum number of double points is $k-1$ and that there are $2^{(k-2)}$ permutations that realize minimal double points, as on our way up we can decide to leave out or visit any of the strands that are not the bottom or the top one (so $k-2$ many yes-no-choices) and on our way down every not yet visited strand has to be visited (no choices). A few more examples for connections with one cycle are given in Figure~\ref{fig:onecycleconnections}

	\begin{figure}[h!]
		\centering
		\includegraphics[scale=0.23]{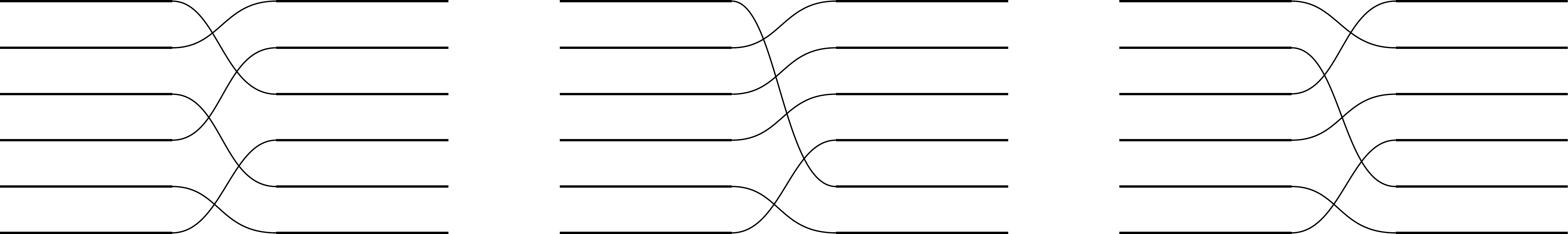}
		\caption{Examples for other connections with one cycle. Here $k = 6 .$}
		\label{fig:onecycleconnections}
	\end{figure}
\end{remark}

Now we are done constructing the $k$-bifurcation for any immersion with the \emph{lowest number of double points}. But we have yet to prove that this construction leads to the lowest $J^+.$ Maybe we could get a lower $J^+$ if we do not connect the strands at the same place or if we use more double points. But we cannot and we show why next.

\begin{thesislemma}
	\label{lem:mindpsminj}
	If $K$ is an immersion with $n_K \in \mathbb{N}$ double points, then a $k$-bifurcation $\widetilde{K}_k$ has minimal $J^+$ if and only if the number of double points is minimal, i.e.\@ the $k$-bifurcation has $ n_K k^2 + (k-1) $ double points.
\end{thesislemma}

The last part of the lemma follows directly from our previous Lemma~\ref{lem:mindps}, from which we know that the minimal number of double points of a $k$-bifurcation of an immersion~$K$ is $n_K k^2 + (k-1) .$ The main part of the lemma follows directly from the next lemma, which states that the value of $J^+$ increases and decreases by the same number as additional double points.

\begin{thesislemma}
	\label{lem:kbifdpsj}
	Let $K$ be an arbitrary immersion and $\widetilde{K}_k, \widetilde{K}'_k$ two arbitrary $k$-bifurcations of $K .$ Then
	$$ J^+(\widetilde{K}'_k) - J^+(\widetilde{K}_k) = n_{\widetilde{K}'_k} - n_{\widetilde{K}_k} .$$
\end{thesislemma}

\begin{proof}[Proof of Lemma~\ref{lem:kbifdpsj}]
	The only double points that can be different between two $k$-bifurcations of $K$ are ones that are between the strands of the bifurcation, which follows from our observations on additional crossings within the path, see Chapter~\ref{subsubsec:addcross}.
	
	We also know that two different $k$-bifurcations are regularly homotopic. We can simply take a linear homotopy between the two bifurcations. Further, small perturbations make the immersion during the homotopy generic for all but isolated moments.
	
	There are three events during regular homotopies that can change the topological properties of an immersion:
	\begin{itemize}
		\item direct self-tangency (see Figure~\ref{fig:jselfdirect})
		\item inverse self-tangency (see Figure~\ref{fig:jselfinverse})
		\item triple point crossing (see Figure~\ref{fig:jtrip})
	\end{itemize}

	A triple point crossing does neither change the number of double points, nor the value of $J^+ .$

	An inverse self-tangency can never happen during a regular homotopy between two $k$-bifurcations of the same immersion~$K ,$ as self-tangencies can only occur between strands of the same part of the path, which are all directed in the same direction, as explored in Chapter~\ref{subsubsec:addcross}.

	So all we need to consider are direct self-tangencies. A positive (negative) direct self-tangency creates (removes) two double points and changes the value of $J^+$ by $2$ ($-2$).
	
	This shows that during a regular homotopy between two different $k$-bifurcations of the same immersion~$K$ the value of $J^+$ and the number of double points change by the same number.\qedhere
\end{proof}

With this we see that the construction from the proof of Lemma~\ref{lem:mindps} in fact has the lowest~$J^+ ,$ as it has minimal number of double points.

\subsection{\texorpdfstring{$J^+$ under $k$-bifurcations}{J+ under k-bifurcations}}
\label{subsec:kbifj}

Now we are ready to prove our main results. First we prove Theorem A, an attainable lower bound for $J^+$ for $k$-bifurcations of any immersion~$K.$ From this Theorem B follows. And then Theorem C, a simple way to calculate $J^+$ of any $k$-bifurcation where only $J^+$ of the original immersion~$K$ and the number of double points of $K$ and its $k$-bifurcation is needed, using the lower bound.

\begin{theooo}[Lower bound for $J^+$ of $k$-bifurcations (Theorem A)]
	\label{th:kbifjlowerbound}
	For all $k$-bifurcations~$\widetilde{K}_k$ of any immersion~$K$ the following lower bound holds:
	\begin{equation*}
		J^+(\widetilde{K}_k) \ge k^2 J^+(K) - (k^2 - k)
	\end{equation*}

	It is attained by any $k$-bifurcation of $K$ with minimal double points.
\end{theooo}

In order to prove this, we use Viro's formula (see Lemma~\ref{lem:viro})
$$J^+(K) = 1 + n - \sum \limits_{C \in \Gamma_K} (\omega_C(K))^2 + \sum \limits_{p \in \mathcal{D}_K} (\operatorname{ind}_p(K))^2$$
and the construction of $k$-bifurcations with minimal $J^+$ that we just explored in the previous chapter: draw every strand to the left (or right) of the previous one and connect them together neatly at one place as in Figure~\ref{fig:neatconnections}.

We denote a $k$-bifurcation of an immersion~$K$ with minimal $J^+$ as $K^*_k$ from here on with arbitrary $k \ge 2 \in \mathbb{N}.$ So what we will do is show that
$$ J^+(K^*_k) = k^2 J^+(K) - (k^2 - k) .$$

Let us use the same notation as in the introduction of Viro's formula for an immersion~$K$ (replace with $K^*_k$ for a $k$-bifurcation with minimal $J^+$):
\begin{itemize}
	\item $n_K$ the number of double points of an immersion~$K$
	\item $\Gamma_K$ the connected components of $\mathbb{C} \setminus K$
	\item $\mathcal{D}_K$ the double points of $K$
	\item $\operatorname{ind}_p(K)$ the index of the double point $p$ in $K ,$ see Definition~\ref{def:indexdp}
\end{itemize}

And for readability, let us introduce the following to rewrite the sums in Viro's formula. The $V$ in the index is to indicate that it is the sum as it is in Viro's formula:
\begin{itemize}
	\item $\Gamma_V(K) \vcentcolon= \sum \limits_{C \in \Gamma_{K}} (\omega_C(K))^2$
	\item $\mathcal{D}_V(K) \vcentcolon= \sum \limits_{p \in \mathcal{D}_{K}} (\operatorname{ind}_p(K))^2$
\end{itemize}

With this, we can write Viro's formula as
$$ J^+(K) = 1 + n_K - \Gamma_V(K) + \mathcal{D}_V(K) .$$

We already know (see Lemma~\ref{lem:mindps}) that
$$n_{K^*_k} = n_K k^2 + (k-1) .$$

To compute $J^+(K^*_k)$ using Viro's formula, we do not compute $\Gamma_V(K^*_k)$ and $\mathcal{D}_V(K^*_k)$ separately. Instead we prove the following lemma and then do a simple calculation.

\begin{thesislemma}
	\label{lem:virosums}
	For all $k$-bifurcations~$K^*_k$ with minimal $J^+$ of any immersion~$K$ the following equation holds:
	\begin{equation*}
		- \Gamma_V(K^*_k) + \mathcal{D}_V(K^*_k) = k^2 \left(- \Gamma_V(K) + \mathcal{D}_V(K) \right)
	\end{equation*}
\end{thesislemma}

\begin{proof}[Proof of Lemma~\ref{lem:virosums}]
	Let us start with an overview of what we will do in this proof before we dive in.

	The idea of this proof is to partition the set of connected components $\Gamma_{K^*_k}$ and the set of double points $\mathcal{D}_{K^*_k}$ into subsets $\Gamma^{\text{old}}_{K^*_k}, \, \Gamma^{\text{square}}_{K^*_k}, \, \Gamma^{\text{stripe}}_{K^*_k}$ and $\mathcal{D}^{\text{old}}_{K^*_k}, \, \mathcal{D}^{\text{square}}_{K^*_k}, \, \mathcal{D}^{\text{stripe}}_{K^*_k},\, \mathcal{D}^{\text{neat}}_{K^*_k}$ respectively. We do this in a way that makes the following easy to show:
	\begin{itemize}
		\item $\Gamma^{\text{old}}_V (K^*_k) = k^2 \cdot \Gamma_V (K)$
		\item $\mathcal{D}^{\text{old}}_V (K^*_k) = k^2 \cdot \mathcal{D}_V (K)$
		\item $\Gamma^{\text{square}}_V (K^*_k) = \mathcal{D}^{\text{square}}_V (K^*_k)$
		\item $\Gamma^{\text{stripe}}_V (K^*_k) = \mathcal{D}^{\text{stripe}}_V (K^*_k) + \mathcal{D}^{\text{neat}}_V (K^*_k)$
	\end{itemize}

	Which already proves our lemma. With the general idea in mind, we now start the proof.

	First we separate \emph{all the connected components $\Gamma_{K^*_k}$} into three categories. Figure~\ref{fig:comps} illustrates them in a $3$-bifurcation of a simple immersion:
	\begin{itemize}
		\item the quadrangles within the $k \cdot k$-square grid of the former double points of the immersion~$K,$ as earlier seen for a $4$-bifurcation in Figure~\ref{fig:doublepointbif}. Denote these connected components as $\Gamma^{\text{square}}_{K^*_k}.$
		\item the stripes between the $k$ many strands. Denote these connected components as $\Gamma^{\text{stripe}}_{K^*_k}.$
		\item all the other connected components. Denote these as $\Gamma^{\text{old}}_{K^*_k}.$
	\end{itemize}

	\begin{figure}[h!]
		\centering
		\includegraphics[scale=0.55]{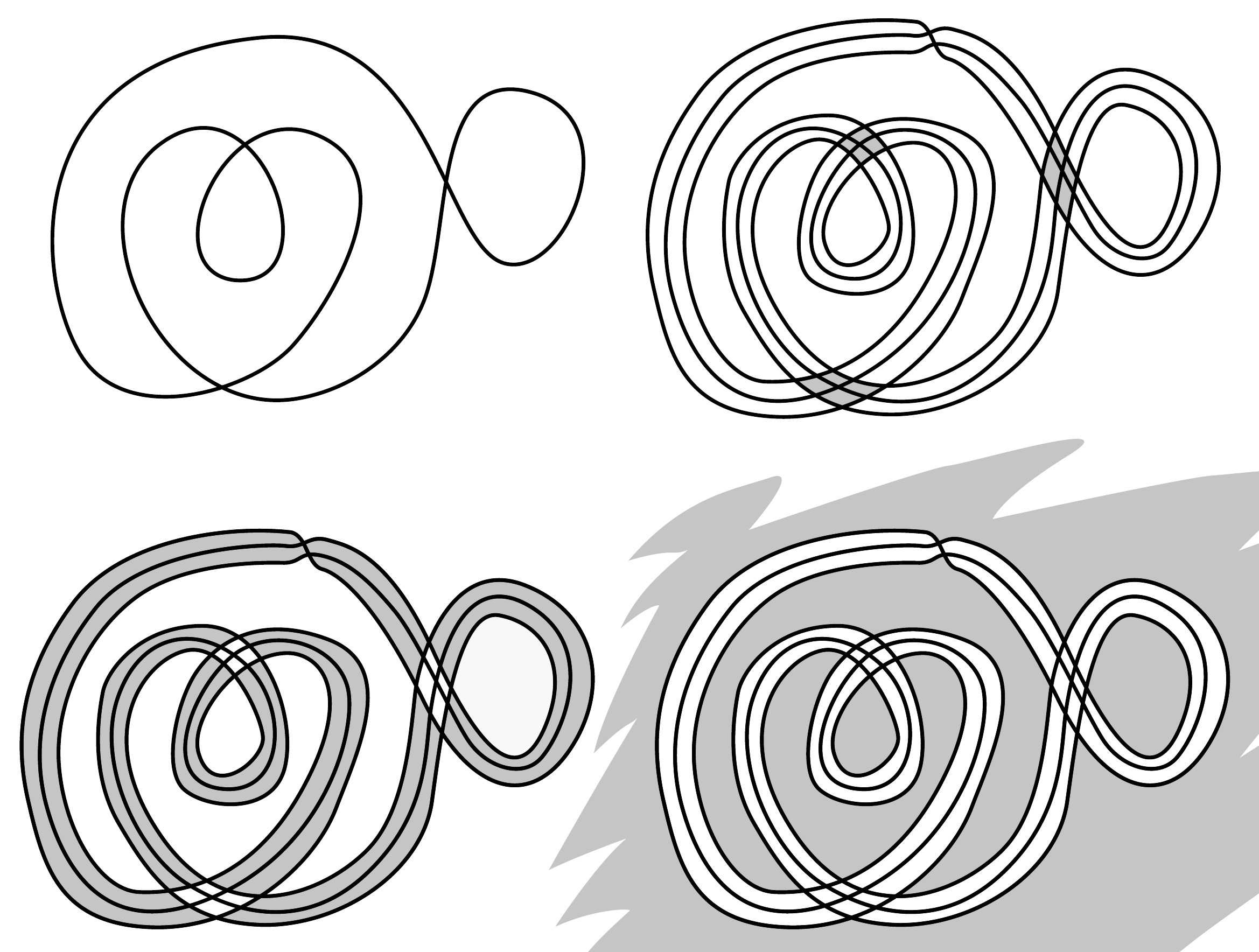}
		\caption{An immersion~$K$ and three identical pictures of a $3$-bifurcation~$K^*_k$ with minimal $J^+.$ Top left: $K.$ Top right: $\Gamma^{\text{square}}_{K^*_k}.$ Bottom left: $\Gamma^{\text{stripe}}_{K^*_k}.$ Bottom right: $\Gamma^{\text{old}}_{K^*_k}.$ Tinted respectively.}
		\label{fig:comps}
	\end{figure}

	The components of $\Gamma^{\text{old}}_{K^*_k}$ are roughly the same area as the connected components~$\Gamma_K$ of $K$ and can be interpreted as the connected components that exist in $\mathbb{C}$ without the path.

	But their winding number, compared to the connected components in $\Gamma_K,$ each is multiplied by~$k.$ This is easy to see using the trick from Chapter~\ref{subsec:imbasics}, Figure~\ref{fig:windingstep}: If a component~$C \in \Gamma_K$ had winding number~$m$, then a flatlander who lives in the plane and starts in the unbounded component with winding number~$0$ would have to cross the immersion transversally $m_+ + m_- \ge m$ times, with~$m_+$~($m_-$) the amount of crossings where the immersion below her went to the right (left) from her point of view. Then the winding number of the component she travelled to is $m = m_+ - m_-.$ Now for each crossing of the immersion~$K$ the flatlander had to do to reach $C,$ in $K^*_k$ she will have to cross $k$ strands all pointing in the same direction (per crossing of the path), so the winding number of that component in $\Gamma^{\text{old}}_{K^*_k}$ will be $k \cdot m_+ - k \cdot m_- = km .$

	So the sum over only these components in $K^*_k$ is
	\begin{align*}
		\Gamma^{\text{old}}_V(K^*_k)
		&= \sum \limits_{C \in \Gamma^{\text{old}}_{K^*_k}} (\omega_C(K^*_k))^2 \\
		&= \sum \limits_{C \in \Gamma_{K}} (k \cdot \omega_C(K))^2 \\
		&= \sum \limits_{C \in \Gamma_{K}} k^2 \cdot (\omega_C(K))^2 \\
		&= k^2 \cdot \sum \limits_{C \in \Gamma_{K}} (\omega_C(K))^2 \\
		&= k^2 \cdot \Gamma_V(K)
	\end{align*}

	Let us denote the connected components in $\Gamma^{\text{stripe}}_{K^*_k}$ as \emph{stripes}. Each stripe is bounded by two neighboring strands and by zero, one or two more strands that are part of the outer edge of a $k^2$-square grid or the one strand that is crossing all the other strands at the neat strand connections. If we look at the orientation of the two neighboring strands, we can use the notion of where the stripe originates from, i.e.\@ what is at the start of the arcs of the neighboring strands that bound the stripe. Obviously, each stripe has exactly one origin: either the neat strand connections or a $k^2$-square grid. There are exactly~$k - 1$ stripes that originate from the neat strand connections and $2 (k - 1)$ stripes that originate from each $k^2$-square grid, see Figure~\ref{fig:stripes}.

	\begin{figure}[h!]
		\centering
		\includegraphics[scale=0.4]{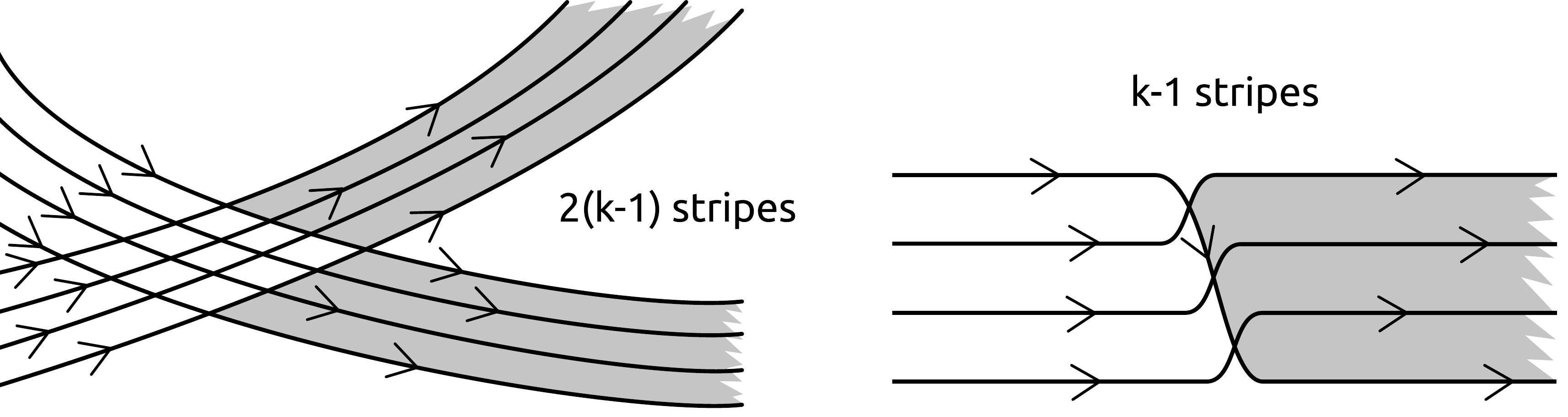}
		\caption{Without loss of generality, we can look at all $k^2$-square grids and at the neat connections like in this figure, with all strands directed from left to right.}
		\label{fig:stripes}
	\end{figure}

	Now let us look at \emph{all the double points $\mathcal{D}_{K^*_k}$.} Here the categories will seem more arbitrary than for the connected components, but they are chosen in a way that will make the rest of the proof easier. Figure~\ref{fig:dps} illustrates them from left to right:

	\begin{itemize}
		\item all of the right most double points of all the $k \cdot k$-square grids when a small neighborhood of the square grid is drawn with all the strands going from left to right, as seen in Figure~\ref{fig:doublepointbif}, i.e.\@ the unique corner of the square grid where the two intersecting strands are oriented away from the square grid. Denote these double points as $\mathcal{D}^{\text{old}}_{K^*_k}.$
		\item all of the double points on the upper right and lower right side of the square grids, with all strands drawn from left to right as before. Exclude the right corner, as it is already in $\mathcal{D}^{\text{old}}_{K^*_k}.$ In other words: all outer double points of the square grids where exactly one of the two intersecting strands is oriented away from the square grid. Denote these double points as $\mathcal{D}^{\text{stripe}}_{K^*_k},$ as they are at the start of the stripes.
		\item all the other double points of the square grids. Denote these double points as $\mathcal{D}^{\text{square}}_{K^*_k}.$
		\item the $k-1$ many double points in the neat strand connections of the $k$ strands as seen in Figure~\ref{fig:neatconnections}. Denote these double points as $\mathcal{D}^{\text{neat}}_{K^*_k}.$
	\end{itemize}

	\begin{figure}[h!]
		\centering
		\includegraphics[scale=0.35]{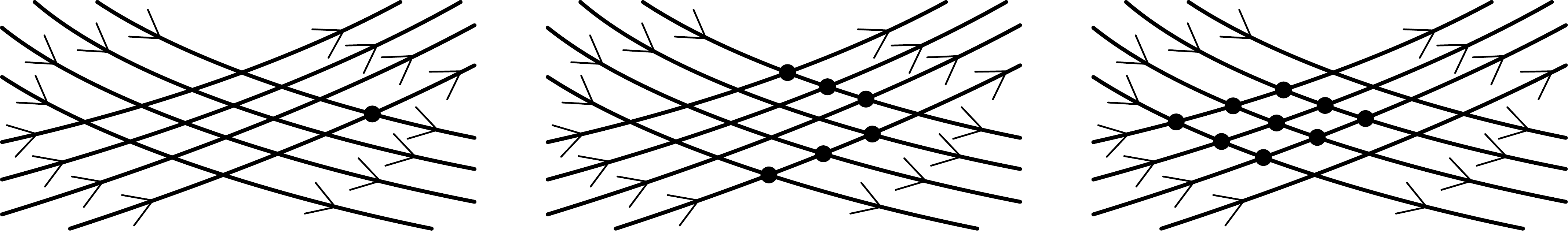}
		\caption{Left to right: $\mathcal{D}^{\text{old}}_{K^*_k}, \, \mathcal{D}^{\text{stripe}}_{K^*_k}, \, \mathcal{D}^{\text{square}}_{K^*_k}$ respectively in a $k^2$-square grid. Here $k=4$ to keep the figure overseeable.}
		\label{fig:dps}
	\end{figure}

	Remember that the index of a double point~$p$ in the pictures we draw -- with all strands directed from left to right -- is equal to the winding number of the connected component to its side, i.e.\@ to its right or left. This follows from our observations of Remark~\ref{rem:indexdpside}.

	The index of a double point in $\mathcal{D}^{\text{old}}_{K^*_k}$ is equal to the winding number of the connected component to the right (or left) of the $k^2$-square grid of the double point, which is equal to $k$ times the index of the double point (from the pre-bifurcation immersion) from which the $k^2$-square grid originates. So similar to the equation for $\Gamma^{\text{old}}_{K^*_k} ,$ for $\mathcal{D}^{\text{old}}_{K^*_k}$ we get:

	\begin{align*}
		\mathcal{D}^{\text{old}}_V(K^*_k)
		&= \sum \limits_{p \in \mathcal{D}^{\text{old}}_{K^*_k}} (\operatorname{ind}_p(K^*_k))^2 \\
		&= \sum \limits_{p \in \mathcal{D}_{K}} (k \cdot \operatorname{ind}_p(K))^2 \\
		&= \sum \limits_{p \in \mathcal{D}_{K}} k^2 \cdot (\operatorname{ind}_p(K))^2 \\
		&= k^2 \cdot \sum \limits_{p \in \mathcal{D}_{K}} (\operatorname{ind}_p(K))^2 \\
		&= k^2 \cdot \mathcal{D}_V(K)
	\end{align*}

	Now we can easily see, as illustrated in Figure~\ref{fig:matchab}, that:

	\begin{itemize}
		\item $\Gamma^{\text{square}}_V (K^*_k) = \mathcal{D}^{\text{square}}_V (K^*_k),$ see the left picture of Figure~\ref{fig:matchab}.
		\item $\Gamma^{\text{stripe}}_V (K^*_k) = \mathcal{D}^{\text{stripe}}_V (K^*_k) + \mathcal{D}^{\text{neat}}_V (K^*_k),$ see the middle and right picture of Figure~\ref{fig:matchab}.
	\end{itemize}

	\begin{figure}[h!]
		\centering
		\includegraphics[scale=0.3]{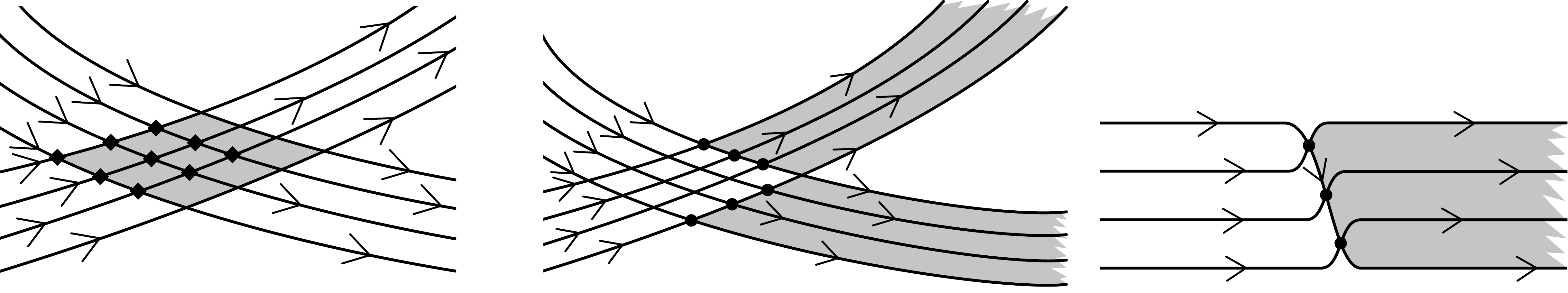}
		\caption{Here it is easy to see which points' index is equal to which connected components' winding number.}
		\label{fig:matchab}
	\end{figure}

	And with this, we showed the statements from the overview at the beginning of this proof. We now use them to finish the proof with this simple calculation:
	\begin{align*}
		- \Gamma_V(K^*_k) + \mathcal{D}_V(K^*_k)
		&= - (\Gamma^{\text{old}}_V(K^*_k) + \Gamma^{\text{square}}_V(K^*_k) + \Gamma^{\text{stripe}}_V(K^*_k)) \\
			&\quad \, + (\mathcal{D}^{\text{old}}_V(K^*_k) + \mathcal{D}^{\text{square}}_V(K^*_k) + \mathcal{D}^{\text{stripe}}_V(K^*_k) + \mathcal{D}^{\text{neat}}_V(K^*_k)) \\
		&= - \Gamma^{\text{old}}_V(K^*_k) + \mathcal{D}^{\text{old}}_V(K^*_k) \\
		&= - k^2 \cdot \Gamma_V(K) + k^2 \cdot \mathcal{D}_V(K) \\
		&= k^2 \left(- \Gamma_V(K) + \mathcal{D}_V(K) \right)
		\qedhere
	\end{align*}
\end{proof}

Now we can use the lemma to prove Theorem A.

\begin{proof}[Proof of Theorem~\ref{th:kbifjlowerbound}: Lower bound for $J^+$ of $k$-bifurcations (Theorem A)]
	\label{pr:kbifjlowerbound}
	Let $K^*_k$ be any $k$-bifurcation of an immersion~$K$ with minimal $J^+,$ so $J^+(\widetilde{K}_k) \ge J^+(K^*_k).$ We show that
	$$ J^+(K^*_k) = k^2 J^+(K) - (k^2 - k) , \forall \widetilde{K}_k, k \ge 2 .$$

	Using Lemma~\ref{lem:virosums} we can calculate
	\begin{align*}
		J^+(K^*_k)
		&= 1 + \underbrace{n_{K^*_k}}_{= n_K k^2 + (k-1)} \underbrace{- \Gamma_V(K^*_k) + \mathcal{D}_V(K^*_k)}_{= k^2 \left(- \Gamma_V(K) + \mathcal{D}_V(K) \right)} \\
		&= 1 + n_K k^2 + (k-1) + k^2 \left(- \Gamma_V(K) + \mathcal{D}_V(K) \right) \\
		&= k^2 \left( 1 + n_K - \Gamma_V(K) + \mathcal{D}_V(K) \right) - (k^2 - k) \\
		&= k^2 J^+(K) - (k^2 - k)
		\qedhere
	\end{align*}
\end{proof}

Let us look at two easy corollaries that follow from the proof of Lemma~\ref{lem:virosums}. We do not need them later, so feel free to skip them.

\begin{thesiscorollary}
	\label{cor:kbifjuselessremark}
	For any immersion~$K$ and $k \ge 2$ we have:
	$$ \Gamma_V(K^*_k) \gneq  k^2 \cdot \Gamma_V(K) $$
	and
	$$ \mathcal{D}_V(K^*_k) \gneq k^2 \cdot \mathcal{D}_V(K) .$$
\end{thesiscorollary}

\begin{proof}
	The winding number of any stripe in $\Gamma^{\text{stripe}}_{K^*_k}$ is $\neq 0.$ Again we use the trick from Figure~\ref{fig:windingstep}. Any stripe is between two connected components of $\Gamma^{\text{old}}_{K^*_k} .$ As we observed at the beginning of the previous proof (of Lemma~\ref{lem:virosums}), any such two components can be interpreted as two formerly adjacent components in~$K, $ so now their winding numbers are two subsequent multiples of $k.$ All strands between these two components must have winding numbers between these subsequent multiples of $k, $ which therefore cannot include $0, $ as it is a multiple of $k$ itself.

	The first half of the corollary follows from this observation.

	The second half then easily follows from this statement from the previous proof:
	\begin{equation*}
		\Gamma^{\text{stripe}}_V (K^*_k) = \mathcal{D}^{\text{stripe}}_V (K^*_k) + \mathcal{D}^{\text{neat}}_V (K^*_k) .
		\qedhere
	\end{equation*}
\end{proof}

\begin{thesiscorollary}
	A few quantitative statements that follow from the proof of Lemma~\ref{lem:virosums}:
	\begin{itemize}
		\item $| \Gamma^{\text{square}}_{K^*_k} | = | \mathcal{D}^{\text{square}}_{K^*_k} | = n_K (k-1)^2$
		\item $| \Gamma^{\text{stripe}}_{K^*_k} | = | \mathcal{D}^{\text{stripe}}_{K^*_k} | + | \mathcal{D}^{\text{neat}}_{K^*_k} | = (2n_K + 1) (k-1)$
	\end{itemize}
\end{thesiscorollary}

With Theorem~\ref{th:kbifjlowerbound} proven, we can now easily conclude:

\begin{thesiscorollary}[$J^+$ stays positive under bifurcations (Theorem B)]
	\label{cor:kbifjpositivek2}
	For all $k$-bifurcations~$\widetilde{K}_k$ of any immersion~$K$ the following holds:
	$$ J^+(K) > 0 \quad \Rightarrow \quad J^+(\widetilde{K}_k) \ge k^2 + k $$
\end{thesiscorollary}

From this it follows that whenever we have one immersion $K_+$ with $J^+(K_+) > 0$ and one immersion $K_-$ with $J^+(K_-) \le 0 ,$ then we know that we cannot get $K_-$ from $K_+$ by applying generic homotopies and bifurcations. Note that $J^+(K) > 0 \Leftrightarrow J^+(K) \ge 2 ,$ because $J^+(K) \in 2 \mathbb{Z}$ for any~$K .$

\begin{proof}[Proof of Corollary~\ref{cor:kbifjpositivek2}: $J^+$ stays positive under bifurcations]
	\label{pr:kbifjpositivek2}
	Let $\widetilde{K}_k$ be any arbitrary $k$-bifurcation ($k \ge 2$) of an arbitrary immersion~$K$ with $J^+(K) > 0 .$

	A simple calculation finishes the proof, using Theorem~\ref{th:kbifjlowerbound}:
	\begin{align*}
		J^+(\widetilde{K}_k) &\ge k^2 J^+(K) - (k^2 - k) \\
		\Leftrightarrow J^+(\widetilde{K}_k) &\ge k^2 ( \underbrace{J^+(K)}_{\ge 2} - 1 ) + k \\
		\Rightarrow J^+(\widetilde{K}_k) &\ge k^2 + k
		\qedhere
	\end{align*}
\end{proof}

Before we observe what follows for the other three invariants $J^- ,$ $\mathcal{J}_1$ and $\mathcal{J}_2 ,$ we take note of the following Theorem, which follows directly from Theorem~\ref{th:kbifjlowerbound} and the preparatory work from Lemma~\ref{lem:kbifdpsj}.

\begin{theooo}[Easy calculation of $J^+$ of any $k$-bifurcation (Theorem C)]
	\label{th:kbifjcalc}
	Let $K$ be an arbitrary immersion and $\widetilde{K}_k$ an arbitrary $k$-bifurcation of $K .$ Then
	$$ J^+(\widetilde{K}_k) = J^+(K^*_k) + (n_{\widetilde{K}_k} - n_{K^*_k}) ,$$
	with $K^*_k$ any $k$-bifurcation of $K$ with minimal $J^+ ,$ and $n_K, n_{\widetilde{K}_k}, n_{K^*_k}$ the number of double points of $K, \widetilde{K}_k, K^*_k$ respectively.
\end{theooo}

So $J^+$ of any $k$-bifurcation of $K$ is equal to $J^+$ of $K^*_k$ -- a $k$-bifurcation of $K$ with minimal number of double points -- plus the number of double points that it has more than $K^*_k .$

\begin{remark}
	Note that
	\begin{align*}
		J^+(K^*_k) + (n_{\widetilde{K}_k} - n_{K^*_k})
		&= k^2 J^+(K) - (k^2 - k) + (n_{\widetilde{K}_k} - (n_K k^2 + (k-1))) \\
		&= k^2 J^+(K) - (k^2 - 1) + n_{\widetilde{K}_k} - n_K k^2
	\end{align*}
\end{remark}

We can now use these findings to formulate results for $J^-$, $\mathcal{J}_1$ and $\mathcal{J}_2$ which we introduced in Chapters~\ref{subsec:jplusminusinvs} and~\ref{subsec:szhom} but have since not used. Their definition is briefly recalled in the following subchapters.

\subsection{\texorpdfstring{$J^-$ under $k$-bifurcations}{J- under k-bifurcations}}
\label{subsec:bifjminus}

Recall that to calculate $J^-$ of an immersion~$K$ we can simply calculate $J^+$ and then subtract the number of double points $n_K ,$ so:
$$ J^+(K) - J^-(K) = n_K . $$

We also know that the additional crossings discussed in Chapter~\ref{subsubsec:addcross} do not change $J^- ,$ as it does not change under direct self-tangencies, instead only under inverse self-tangencies. With this in mind we see that for any arbitrary immersion $K$ and two arbitrary $k$-bifurcations $\widetilde{K}_k, \widetilde{K}'_k$ of $K,$ we have
$$ J^-(\widetilde{K}'_k) - J^-(\widetilde{K}_k) = 0 ,$$
unlike for $J^+$ (see Lemma~\ref{lem:kbifdpsj}).

This means that $J^-$ is the same for any $k$-bifurcation (with fixed $k$). With this info we can formulate the following Corollary~\ref{cor:kbifjminus}, either by using Theorem~\ref{th:kbifjcalc} or with the following proof.

\begin{thesiscorollary}[$J^-$ of $k$-bifurcations]
	\label{cor:kbifjminus}
	For all $k$-bifurcations~$\widetilde{K}_k$ of any immersion~$K$ we have:
	\begin{equation*}
		J^-(\widetilde{K}_k) = k^2 J^-(K) - (k^2 - 1)
	\end{equation*}
\end{thesiscorollary}

The proof is similar to the proof for the lower bound of $J^+ ,$ also using Lemma~\ref{lem:virosums}.

\begin{proof}
	Let $\widetilde{K}_k$ be any $k$-bifurcation of an immersion~$K .$

	Using Lemma~\ref{lem:virosums} we can calculate
	\begin{align*}
		J^-(\widetilde{K}_k)
		&= 1 \underbrace{- \Gamma_V(\widetilde{K}_k) + \mathcal{D}_V(\widetilde{K}_k)}_{= k^2 \left(- \Gamma_V(K) + \mathcal{D}_V(K) \right)} \\
		&= 1 + k^2 \left(- \Gamma_V(K) + \mathcal{D}_V(K) \right) \\
		&= k^2 \left( 1 - \Gamma_V(K) + \mathcal{D}_V(K) \right) - k^2 + 1 \\
		&= k^2 J^-(K) - (k^2 - 1)
		\qedhere
	\end{align*}
\end{proof}

A quick calculation gives us the following result.

\begin{thesiscorollary}[$J^-$ stays positive (for positive integers) under bifurcations]
	For all $k$-bifurcations~$\widetilde{K}_k$ of any immersion~$K$ the following holds:
	$$ J^-(K) \ge 1 \quad \Rightarrow \quad J^-(\widetilde{K}_k) \ge 1 $$
\end{thesiscorollary}

\begin{proof}
	Using Corollary~\ref{cor:kbifjminus}:
	\begin{align*}
		J^-(\widetilde{K}_k)
		&= k^2 J^-(K) - (k^2 - 1) \\
		&= k^2 (J^-(K) - 1) + 1 \\
		&\ge k^2 (1 - 1) + 1 \\
		&= 1
		\qedhere
	\end{align*}
\end{proof}

\subsection{\texorpdfstring{$\mathcal{J}_1$ under $k$-bifurcations}{J1 under k-bifurcations}}
\label{subsec:bifj1}

Recall that to calculate $\mathcal{J}_1$ of an immersion~$K$ there needs to be a specified origin point~$0 \in \mathbb{C} \setminus K$ and~$K \subset \mathbb{C}^* .$ Then we defined
$$\mathcal{J}_1(K) \vcentcolon= J^+(K) + \dfrac{\omega_0(K)^2}{2} ,$$
see Definition~\ref{def:j1}.

Using additional crossings, we can easily see that $\mathcal{J}_1$ has no upper bound, just like $J^+.$ And just like for $J^+$ we can find a lower bound:

\begin{thesiscorollary}[Lower bound for $\mathcal{J}_1$ of $k$-bifurcations]
	\label{cor:kbifj1lowerbound}
	For all $k$-bifurcations~$\widetilde{K}_k$ of any immersion~$K$ the following lower bound holds:
	\begin{equation}
		\label{eq:kbifj1lowerbound}
		\mathcal{J}_1(\widetilde{K}_k) \ge k^2 \mathcal{J}_1(K) - (k^2 - k)
	\end{equation}

	It is attained by any $k$-bifurcation of $K$ with minimal double points.
\end{thesiscorollary}

\begin{proof}
	\label{pr:kbifj1lowerbound}
	This follows from Theorem~\ref{th:kbifjlowerbound} and a simple calculation.

	Like in Chapter~\ref{subsec:kbifj}, denote a $k$-bifurcation of $K$ with minimal $J^+$ as $K^*_k .$ From the proof of Lemma~\ref{lem:virosums} we know that $\omega_0(K^*_k) = \omega_0(\widetilde{K}_k) = k \cdot \omega_0(K) .$ This short calculation finishes the proof:
	\begin{align*}
		\mathcal{J}_1(K^*_k)
		&= \underbrace{J^+(K^*_k)}_{= k^2 \cdot J^+(K) - (k^2 - k)} + \dfrac{\omega_0(K^*_k)^2}{2} \\
		&= (k^2 \cdot J^+(K) + \dfrac{k^2 \cdot \omega_0(K)^2}{2}) - (k^2 - k) \\
		&= k^2(J^+(K) + \dfrac{\omega_0(K)^2}{2}) - (k^2 - k) \\
		&= k^2(\mathcal{J}_1(K)) - (k^2 - k)
		\qedhere
	\end{align*}
\end{proof}

Further, from Lemma~\ref{lem:kbifdpsj} and Theorem~\ref{th:kbifjcalc} we also get the following corollary.

\begin{thesiscorollary}[Easy calculation of $\mathcal{J}_1$ of any $k$-bifurcation]
	\label{cor:kbifj1calc}
	Let $K$ be an arbitrary immersion and $\widetilde{K}_k$ an arbitrary $k$-bifurcation of $K .$ Then
	$$\mathcal{J}_1(\widetilde{K}_k) = \mathcal{J}_1(K^*_k) + (n_{\widetilde{K}_k} - n_{K^*_k}) ,$$
	with $K^*_k$ any $k$-bifurcation of $K$ with minimal $\mathcal{J}_1 ,$ and $n_K, n_{\widetilde{K}_k}, n_{K^*_k}$ the number of double points of $K, \widetilde{K}_k, K^*_k$ respectively.
\end{thesiscorollary}
\vspace*{0.8em}

So $\mathcal{J}_1$ of any $k$-bifurcation of $K$ is equal to $\mathcal{J}_1$ of $K^*_k$ -- a $k$-bifurcation of $K$ with minimal number of double points -- plus the number of double points that it has more than $K^*_k .$

\begin{remark}
	Note that
	\begin{align*}
		\mathcal{J}_1(K^*_k) + (n_{\widetilde{K}_k} - n_{K^*_k})
		&= k^2 \mathcal{J}_1(K) - (k^2 - k) + (n_{\widetilde{K}_k} - (n_K k^2 + (k-1))) \\
		&= k^2 \mathcal{J}_1(K) - (k^2 - 1) + n_{\widetilde{K}_k} - n_K k^2
	\end{align*}
\end{remark}

And of course we can also get the following corollary from Corollary~\ref{cor:kbifj1lowerbound} similar to how we got Corollary~\ref{cor:kbifjpositivek2} from Theorem~\ref{th:kbifjlowerbound}.

\begin{thesiscorollary}[$\mathcal{J}_1$ stays positive (for positive integers) under bifurcations]
	\label{cor:kbifj1positive}
	For all $k$-bifurcations~$\widetilde{K}_k$ of any immersion~$K$ the following holds:
	$$ \mathcal{J}_1(K) > \frac{1}{2} \quad \Rightarrow \quad \mathcal{J}_1(\widetilde{K}_k) \ge k^2 + k $$
\end{thesiscorollary}

Note that $\mathcal{J}_1$ can be $\frac{1}{2}$ (e.g.\@ for the standard curve $K_3$ with the origin in one of the loops), but the corollary only works for $\mathcal{J}_1(K) \ge 2 .$

\subsection{\texorpdfstring{$\mathcal{J}_2$ under $k$-bifurcations}{J2 under k-bifurcations}}
\label{subsec:bifj2}

Recall that to calculate $\mathcal{J}_2$ of an immersion~$K$ there needs to be a specified origin point~$0 \in \mathbb{C} \setminus K$ and~$K \subset \mathbb{C}^* .$ Then we defined
$$\mathcal{J}_2(K) \vcentcolon= \begin{cases}
	J^+(L^{-1} (K)), &\text{if \( \omega_0(K) \) odd} \\
	J^+(\widehat{K}), &\text{else}
\end{cases} ,$$
with $L: \mathbb{C}^* \to \mathbb{C}^*, v \mapsto v^2$ the complex squaring map and $\widehat{K}$ one of the two up to rotation identical lifts of $K ,$ see Chapter~\ref{subsubsec:liftsdpsj2}.

Remember that $J^+$ is not simply a summand in the calculation of $\mathcal{J}_2 .$ Yet we only need a few observations to prove similar lower bounds and easy calculation for $\mathcal{J}_2$ of $k$-bifurcations. It \emph{is} a little different from $J^+$ and $\mathcal{J}_1 ,$ but not too much. This time we get different equations depending on whether the winding number of the immersion~$K$ is even and the winding number of its $k$-bifurcation~$\widetilde{K}_k$ is even -- the latter follows if $\omega_0(K)$ or $k$ is even.

\begin{propooo}
	\label{prop:kbif2lowerbound}
	For all $k$-bifurcations~$\widetilde{K}_k$ of any immersion~$K$ the following lower bounds hold depending on the parity of $\omega_0(K)$ and $k$:
	$$ \mathcal{J}_2(\widetilde{K}_k) \ge \begin{cases}
		k^2 \mathcal{J}_2(K) - (k^2 - k) , &\text{if \( \omega_0(K) \) even} \\
		k^2 \mathcal{J}_2(K) - (k^2 - k) + (k - 1) , &\text{if \( \omega_0(K) \) odd, \( k \) odd} \\
		(\frac{k}{2})^2 \mathcal{J}_2(K) - ((\frac{k}{2})^2 - \frac{k}{2}) , &\text{if \( \omega_0(K) \) odd, \( k \) even}
	\end{cases} $$

	The lower bound for if \emph{$\omega_0(K)$ is odd and $k$ is even} is attained by any $k$-bifurcation with minimal number of double points in a lift of $K$ -- which is equivalent to having minimal number of even double points in $K ,$ see Definition~\ref{def:evendp}.

	The lower bounds for the other two cases are attained only by any $k$-bifurcation with minimal number of double points.
\end{propooo}
\vspace*{0.8em}

We prove this together with the following proposition. Similarly to $J^+$ and $\mathcal{J}_1 ,$ with our results from Lemma~\ref{lem:kbifdpsj} and Theorem~\ref{th:kbifjcalc}, we can also prove a way to calculate $\mathcal{J}_2$ of any $k$-bifurcation.

\begin{propooo}[Easy calculation of $\mathcal{J}_2$ of any $k$-bifurcation]
	\label{prop:kbif2calc}
	Let $K$ be an arbitrary immersion and $\widetilde{K}_k$ an arbitrary $k$-bifurcation of $K .$ Then
	$$ \mathcal{J}_2(\widetilde{K}_k) = \begin{cases}
		\mathcal{J}_2(K^*_k) + (n_{\widetilde{K}_k} - n_{K^*_k}) , &\text{if \( \omega_0(K) \) even} \\
		\mathcal{J}_2(K^*_k) + 2(n_{\widetilde{K}_k} - n_{K^*_k}) , &\text{if \( \omega_0(K) \) odd, \( k \) odd} \\
		\mathcal{J}_2(K^*_k) + (\nu_{\widetilde{K}_k} - \nu_{K^*_k}) , &\text{if \( \omega_0(K) \) odd, \( k \) even}
	\end{cases} $$
	with $K^*_k$ any $k$-bifurcation of $K$ with minimal $\mathcal{J}_2 ,$ and $n_K, n_{\widetilde{K}_k}, n_{K^*_k}$ the number of double points of $K, \widetilde{K}_k, K^*_k$ respectively. We denote the number of even double points (see Definition~\ref{def:evendp}) of an immersion~$K$ with winding number $\omega_0(K)$ even as $\nu_K ,$ which is equal to the number of double points of any of its two lifts.
\end{propooo}

\begin{remark}
	Note that for $\omega_0(K)$ even:
	\begin{align*}
		\mathcal{J}_2(K^*_k) + (n_{\widetilde{K}_k} - n_{K^*_k})
		&= k^2 \mathcal{J}_2(K) - (k^2 - k) + (n_{\widetilde{K}_k} - (n_K k^2 + (k-1))) \\
		&= k^2 \mathcal{J}_2(K) - (k^2 - 1) + n_{\widetilde{K}_k} - n_K k^2 .
	\end{align*}

	For $\omega_0(K)$ odd, $k$ odd:
	\begin{align*}
		\mathcal{J}_2(K^*_k) + 2(n_{\widetilde{K}_k} - n_{K^*_k})
		&= k^2 \mathcal{J}_2(K) - (k^2 - k) + (k - 1) + 2(n_{\widetilde{K}_k} - (n_K k^2 + (k-1))) \\
		&= k^2 \mathcal{J}_2(K) - (k^2 - 1) + 2 (n_{\widetilde{K}_k} - n_K k^2) .
	\end{align*}

	For $\omega_0(K)$ odd, $k$ even:
	\begin{align*}
		\mathcal{J}_2(K^*_k) + (\nu_{\widetilde{K}_k} - \nu_{K^{*}_k})
		&= \mathcal{J}_2(K^{**}_k) + (\nu_{\widetilde{K}_k} - \nu_{K^{**}_k}) \\
		&= (\textstyle\frac{k}{2})^2 \mathcal{J}_2(K) - ((\textstyle\frac{k}{2})^2 - \textstyle\frac{k}{2}) + (\nu_{\widetilde{K}_k} - \nu_{K^{**}_k}) ,
	\end{align*}
	with $K^{**}_k$ a $k$-bifurcation of $K$ with minimal number of even double points.

	Obviously a $k$-bifurcation~$K^*_k$ (with minimal number of double points) of an immersion $K$ also has minimal number of even double points,
	$$\nu_{K^{**}_k} = \nu_{K^*_k} ,$$
	so using $K^*_k$ in the calculation of this case is unnecessarily strict, but not wrong.

	And we can show that, for the case $\omega_0(K)$ odd, $k$ even:
	\begin{equation}
		\label{eq:practicallyuselessequationaboutevendoublepointsofbifurcations}
		\nu_{K^{**}_k} = \nu_{K^*_k} = n_K \textstyle\frac{k^2}{2} + (\textstyle\frac{k}{2} - 1)
	\end{equation}
	But this information is probably practically useless when calculating $\mathcal{J}_2$ of bifurcations ($\omega_0(K)$ odd, $k \ge 4$ even) by hand, so we give no complete proof. This information is especially useless if we ever need to really count the number of even double points of a bifurcation by hand and we happen to have one of the lifts visualized, because then we can instead just count the number of all double points of the lift. This is because from Corollary~\ref{cor:evenliftdps} we know that for any immersion $K$ with winding number $\omega_0(K)$ even we have $\nu_K = n_{\widehat{K}} .$

	But if we do not have the lift at hand and can only work on the bifurcation itself, then maybe some efficient algorithm can use this information somehow. Let us hope that someday someone can find some nice trick to reduce the number of double points to check for their parity. For anyone interested though, here is the idea of the proof of Equation~\ref{eq:practicallyuselessequationaboutevendoublepointsofbifurcations}:
	
	Consider the bifurcation~$K^*_k$ with minimal number of double points and the neat connections. At each $k^2$-square grid, we can see that for every double point $p$ the (two, three or four) neighbors all have exactly the other parity. So each $k^2$-square grid looks like either one of the first two pictures in Figure~\ref{fig:oddevennu}.
	
	\begin{figure}[h!]
		\centering
		\includegraphics[scale=0.3]{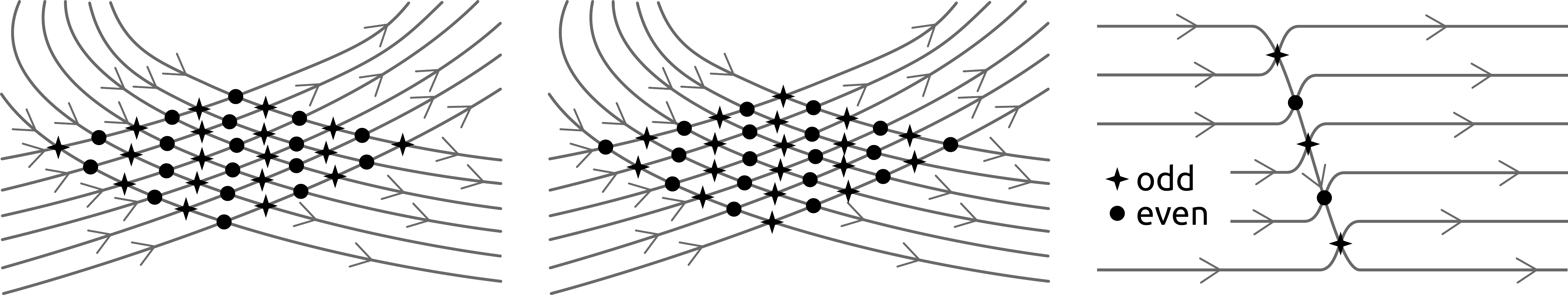}
		\caption{Parity of the double points of a minimal double points $k$-bifurcation~$K^*_k$ for the case $\omega_0(K)$ odd, $k$ even. Here $k = 6$ is visualized.}
		\label{fig:oddevennu}
	\end{figure}
	
	This happens because if we compare the decompositions $K_p^A$ and $K_{p'}^A ,$ (with $p'$ any of the neighboring double points of $p$ and $K_p^A, K_{p'}^A$ both starting from their points in the same direction of the path), we have to pass the path (see Definition~\ref{def:pathstrands}) either once or $(k-1)$ times less or more often in the decomposition of the neighbor. And because $\omega_0(K)$ is odd, this results in a change of parity in either case.
	
	For the neat connections we can see that the same argument (but a little easier) leads to the observation that the double points alternate between odd and even, with odd double points on the outside, like in the right picture of Figure~\ref{fig:oddevennu}. With this, we get that there are $n_K \textstyle\frac{k^2}{2} + (\textstyle\frac{k}{2} - 1)$ even double points for $K^*_k .$
\end{remark}

We prove the three cases of both propositions separately. For all cases let $K$ be an arbitrary immersion,~$\widetilde{K}_k$ an arbitrary $k$-bifurcation of $K$ and $K^*_k$ any $k$-bifurcation of $K$ with minimal $J^+,$ so $J^+(\widetilde{K}_k) \ge J^+(K^*_k) ,$ as always so far.

\begin{proof}[Proof of Propositions~\ref{prop:kbif2lowerbound} and~\ref{prop:kbif2calc}, case $\omega_0(K)$ odd, $k$ odd]
	We start with the easiest case, $\omega_0(K)$ odd, $k$ odd. This is easy to prove with \cite[Proposition 6]{kai:paper}, which states that for any arbitrary immersion $K$ with $\omega_0(K)$ odd the equation
	$$ \mathcal{J}_2(K) = 2 \mathcal{J}_1(K) - 1 $$
	holds. From Corollary~\ref{cor:kbifj1lowerbound} we know that the lower bound of $\mathcal{J}_1(\widetilde{K}_k) ,$ which is $\mathcal{J}_1(\widetilde{K}_k) \ge k^2 \mathcal{J}_1(K) - (k^2 - k)$ (see Equation~\ref{eq:kbifj1lowerbound}), is only attained by $K^*_k .$

	With this we calculate
	\begin{align*}
		\mathcal{J}_2(K^*_k)
		&= 2 \mathcal{J}_1(K^*_k) - 1 \\
		&= 2 (k^2 \mathcal{J}_1(K) - (k^2 - k)) - 1 \\
		&= k^2 (\underbrace{2 \mathcal{J}_1(K) - 1}_{= \mathcal{J}_2(K)} + 1) - 2 (k^2 - k) - 1 \\
		&= k^2 \mathcal{J}_2(K) + k^2 - 2 k^2 + 2 k - 1 \\
		&= k^2 \mathcal{J}_2(K) - (k^2 - k) + (k - 1)
	\end{align*}

	This shows the lower bound for this case. With Corollary~\ref{cor:kbifj1calc} we can calculate the following to show Proposition~\ref{prop:kbif2calc} for this case:
	\begin{align*}
		\mathcal{J}_2(\widetilde{K}_k)
		&= 2 \mathcal{J}_1(\widetilde{K}_k) - 1 \\
		&= 2 (k^2 \mathcal{J}_1(K) - (k^2 - k) + (n_{\widetilde{K}_k} - n_{K^*_k})) - 1 \\
		&= k^2 (\underbrace{2 \mathcal{J}_1(K) - 1}_{= \mathcal{J}_2(K)} + 1) - 2 (k^2 - k) + 2 (n_{\widetilde{K}_k} - n_{K^*_k}) - 1 \\
		&= k^2 \mathcal{J}_2(K) + k^2 - 2 k^2 + 2 k - 1 + 2 (n_{\widetilde{K}_k} - n_{K^*_k}) \\
		&= k^2 \mathcal{J}_2(K) - (k^2 - k) + (k - 1) + 2 (n_{\widetilde{K}_k} - n_{K^*_k}) \\
		&= \mathcal{J}_2(K^*_k) + 2 (n_{\widetilde{K}_k} - n_{K^*_k})
		\qedhere
	\end{align*}
\end{proof}

\begin{remark}
	The case $\omega_0(K)$ odd, $k$ odd can also be proven without \cite[Proposition 6]{kai:paper}. If we recall the proof of Lemma~\ref{lem:virosums} and Theorem~\ref{th:kbifjlowerbound} we see that we need to add $(k-1)$ at the end because the neat connections of $K$ appear twice in the preimage, see our observations from Chapter~\ref{subsubsec:liftsdpsj2}. The difference of double points has to be added twice in the easy calculation for the same reason: each additional crossing in $\widetilde{K}_k$ appears twice in its preimage.
\end{remark}

For the next two cases we denote the lift of an immersion $K$ as $\widehat{K} ,$ like in Chapter~\ref{subsubsec:liftsdpsj2}, and the number of even double points of an immersion~$K$ as $\nu_K .$

\begin{proof}[Proof of Propositions~\ref{prop:kbif2lowerbound} and~\ref{prop:kbif2calc}, case $\omega_0(K)$ even]
	Next we prove the case for $\omega_0(K)$ even, independent of whether $k$ is odd or even.
	
	We observe that for the $k$-bifurcation $K^*_k$ with minimal double points we can calculate $\mathcal{J}_2$ exactly the same way as $J^+ ,$ because the bifurcation of $K$ directly translates to both lifts of $K ,$ which stay separated because the winding number $\omega_0(K^*_k)$ after the bifurcation stays even. This proves the lower bound
	$$ \mathcal{J}_2(\widetilde{K}_k) \ge k^2 \mathcal{J}_2(K) - (k^2 - k) .$$

	For the easy calculation of any $k$-bifurcation $\widetilde{K}_k ,$ we use Theorem~\ref{th:kbifjcalc} on a lift of $\widetilde{K}_k .$ So we need to find the difference of double points between a lift of $K^*_k$ and a lift of $\widetilde{K}_k .$ This means we need to find out which additional crossings in $\widetilde{K}_k$ result in \emph{even} double points, because they are the double points that will also appear in the lift, as we explored earlier in Chapter~\ref{subsubsec:liftsdpsj2}.

	It turns out that for the case $\omega_0(K)$ even, all additional crossings result in even double points. This is because double points from additional crossings are always double points between strands and any of the two decompositions (see Definition~\ref{def:evendp}) at such a double point $p$ is just $j \in \{ 1, \dots, k-1 \}$ runs through the path, which has even winding number.

	This shows that all additional double points in $\widetilde{K}_k$ appear in its lift, so
	$$ (n_{\widehat{\widetilde{K}}_k} - n_{\widehat{K^*_k}}) = (\nu_{\widetilde{K}_k} - \nu_{K^*_k}) = (n_{\widetilde{K}_k} - n_{K^*_k}) $$
	and this shows the easy calculation for this case.
\end{proof}

\begin{figure}[h!]
	\centering
	\includegraphics[scale=0.3]{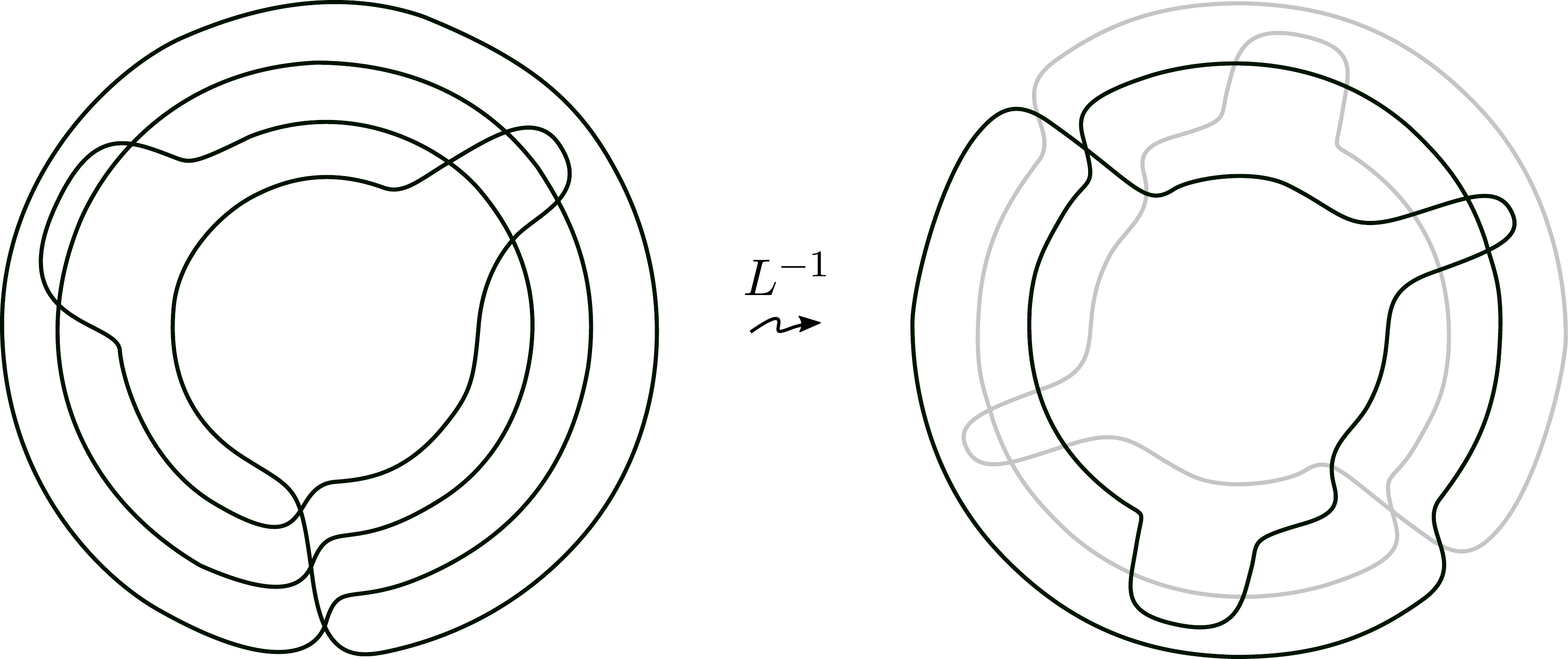}
	\caption{Neat $4$-bifurcation of a circle and its preimage of the squaring map $L$ with some additional crossings. The squaring of the absolute value is ignored here. Note that the figures in this subchapter are hand-drawn, not computer-generated, so not 100\% accurate.}
	\label{fig:evenoddaddcrossings}
\end{figure}

\begin{proof}[Proof of Propositions~\ref{prop:kbif2lowerbound} and~\ref{prop:kbif2calc}, case $\omega_0(K)$ odd, $k$ even]
	The last case is $\omega_0(K)$ odd, $k$ even. The argument here appears simple, and we think it is, but in this paper we lack the proper ideas to perform this proof more precisely.

	The idea is that because the strands of a bifurcation all follow the same path, we basically duplicate the lifts through a $2$-bifurcation, so now we have two closed lifts in $K^*_2$, which each have the same~$J^+$-value as the two combined lifts of $K$ have. For each further increase of $2$ of the bifurcation, the lifts of the bifurcation get one more strand. This brings us to the same form for the lower bound as for $J^+ ,$ but with $\textstyle\frac{k}{2}$ instead of $k .$

	For the easy calculation, like in the previous proof, we just need to add the difference of even double points. This already shows the easy calculation for this case.
	
	But it is important to observe, that other than in the previous proof, the difference of even double points does not turn out to be equal to the difference of \emph{all} double points. So
	$$ (\nu_{\widetilde{K}_k} - \nu_{K^*_k}) \neq (n_{\widetilde{K}_k} - n_{K^*_k}) $$

	This is because the argument of the previous proof does not work here. In general additional crossings do now always create even double points. In fact we can observe with the same argument as in the previous proof, that any additional crossing between neighboring strands of $K^*_k$ creates odd double points. Only additional crossings between strands that have an odd number of other strands between them in $K^*_k$ create even double points. See Figure~\ref{fig:evenoddaddcrossings} for an easy example.
\end{proof}
	
We take note of a single special case of immersions undergoing a bifurcation, where we not only obtain a lower bound, but the invariant $\mathcal{J}_2$ stays exactly the same without any changes to it.

\begin{thesiscorollary}[Invariance of $\mathcal{J}_2$ for $2$-bifurcations of immersions with odd winding number]
	\label{cor:oddj2bif}
	For all $2$-bifurcations~$\widetilde{K}_2$ of any immersion~$K$ with winding number around the origin $\omega_0(K)$ odd, the following equation holds:
	\begin{equation*}
		\mathcal{J}_2(\widetilde{K}_2) = \mathcal{J}_2(K)
	\end{equation*}
\end{thesiscorollary}

So the $\mathcal{J}_2$-invariant of any immersion with odd winding number will stay exactly the same under $2$-bifurcations, no matter how many additional crossings occur.

\begin{proof}
	This follows directly from Proposition~\ref{prop:kbif2calc} and the arguments from the proof of Proposition~\ref{prop:kbif2calc} for this case, which is $\omega_0(K)$ odd, $k$ even.

	With this we get
	\begin{equation*}
		\mathcal{J}_2(\widetilde{K}_k) = \underbrace{\mathcal{J}_2(K^{**}_2)}_{= \mathcal{J}_2(K)} + \underbrace{(\nu_{\widetilde{K}_2} - \nu_{K^{**}_k})}_{= 0} = \mathcal{J}_2(K) .
		\qedhere
	\end{equation*}
\end{proof}

And of course we can also get the following corollary from Proposition~\ref{prop:kbif2lowerbound} just like we got Corollary~\ref{cor:kbifjpositivek2} from Theorem~\ref{th:kbifjlowerbound}.

\begin{thesiscorollary}[$\mathcal{J}_2$ stays positive under bifurcations]
	\label{cor:kbifj2positive}
	For all $k$-bifurcations~$\widetilde{K}_k$ of any immersion~$K$ the 
	following holds:

	$$ \mathcal{J}_2(K) > 0 \quad \Rightarrow \quad \mathcal{J}_2(\widetilde{K}_k) \ge \begin{cases}
		k^2 + k , &\text{if \( \omega_0(K) \) even or \( k \) odd} \\
		(\frac{k}{2})^2 + \frac{k}{2} , &\text{if \( \omega_0(K) \) odd and \( k \) even}
	\end{cases} $$
\end{thesiscorollary}

% \clearpage

\section*{Appendix}
\addcontentsline{toc}{section}{\hspace{1.4em}Appendix}

\fancyhead[RO, LE]{Appendix}
\fancyhead[LO, RE]{}

Welcome to the appendix. The tex files and all image files of this paper can be found in the Git repository on \href{https://gitlab.com/CptMaister/paper-j-bifurcations}{\textbf{gitlab.com/CptMaister/paper-j-bifurcations}} and are available for anyone to use and modify. This paper is based on the bachelor thesis ``Bifurcations of Periodic Orbits in Hamiltonian Systems and $J^+$-Invariants'', which can be found in the Git repository on \href{https://gitlab.com/CptMaister/bachelorthesis}{\textbf{gitlab.com/CptMaister/bachelorthesis}}.

Visit the website \href{http://complex-fibers.org/}{\textbf{complex-fibers.org}} to visualize preimages of the complex squaring map. It can be of great help when trying to understand the $\mathcal{J}_2$-invariant.

\subsection*{Example Viro's Formula}
\addcontentsline{toc}{subsection}{\hspace{2.3em}Example Viro's Formula}

\paragraph*{}

\begin{figure}[h!]
	\centering
	\includegraphics[scale=0.45]{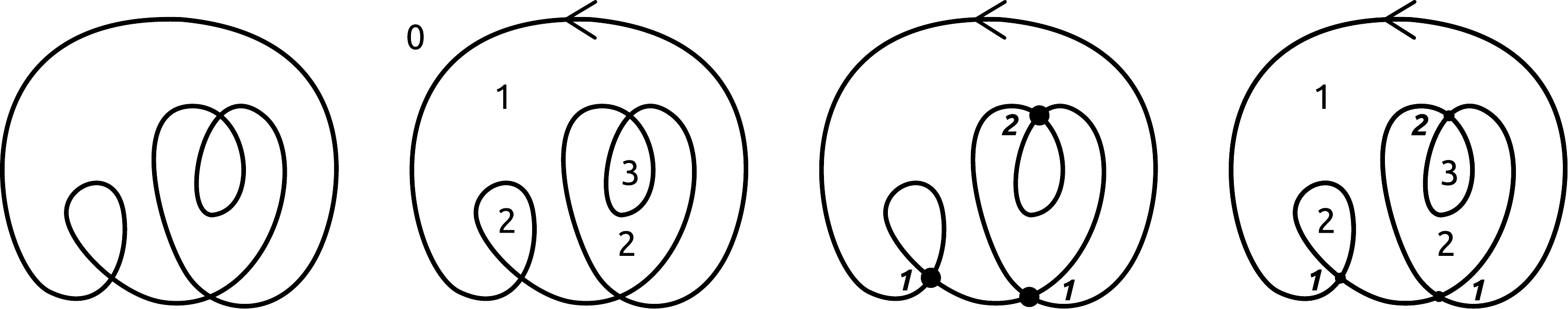}
	\caption{A circle with one single interior loop and one double interior loop. Picture one is just the immersion without orientation, picture two with winding numbers, picture three with double point indices, picture four with both.}
	\label{fig:viroex}
\end{figure}

Let us use Viro's formula to calculate $J^+$ of the immersion~$K$ from Figure~\ref{fig:viroex}. It is a circle with one single interior loop and one double interior loop.\label{example:viro}

In order to use Viro's formula, we need to count the number of double points, get the square of the winding number for all connected components of $\mathbb{C} \setminus K$ and the square of the double point indices.

The number of double points is $3 ,$ which can be easily seen in the third picture of Figure~\ref{fig:viroex}, where the double points are marked. So
$$ n_K = 3 . $$

Next we need all winding numbers. We started with an immersion~$K$ that has no orientation. For Viro's formula we will need the square of all winding numbers of the connected components, so it does not matter which orientation we choose. In this example we choose it so that the immersion has positive rotation number. With the orientation chosen, we can label all winding numbers of the connected components of~$\mathbb{C} \setminus K ,$ see the second picture of Figure~\ref{fig:viroex}. With this we can calculate
\begin{align*}
	\sum \limits_{C \in \Gamma_K} (\omega_C(K))^2
	&= 0^2 + 1^2 + 2^2 + 2^2 + 3^2 \\
	&= 0 + 1 + 4 + 4 + 9 \\
	&= 18 .
\end{align*}

Now the last thing we need is the double point indices. Again we only need the squares of them, so it is ok that we choose any of the two orientations. Remember that the index of a double point is the same as the winding number that appears twice around that double point (see Remark~\ref{rem:indexdptwice}). With this we can easily label all double point indices, see the third picture of Figure~\ref{fig:viroex}. With this we can calculate
\begin{align*}
	\sum \limits_{p \in \mathcal{D}_K} (\operatorname{ind}_p(K))^2
	&= 1^2 + 1^2 + 2^2 \\
	&= 1 + 1 + 4 \\
	&= 6 .
\end{align*}

We combine these observations to get $J^+(K)$ using Viro's formula:
\begin{align*}
	J^+(K)
	&= 1 + \underbrace{n_K}_{= 3} - \underbrace{\sum \limits_{C \in \Gamma_K} (\omega_C(K))^2}_{= 18} + \underbrace{\sum \limits_{p \in \mathcal{D}_K} (\operatorname{ind}_p(K))^2}_{= 6} \\
	&= 1 + 3 - 18 + 6 \\
	&= -8 .
\end{align*}

\clearpage

\pagestyle{litpage}
\nocite{*}
\section*{References}
\addcontentsline{toc}{section}{\hspace{1.4em}References}
\printbibliography[
heading=none,
title={References}
]

\end{document}